%% file: RSC-CIKG-NRL-TN-arXiv.tex
\documentclass[11pt, letterpaper, twoside]{scrartcl}
\usepackage{mathtools}
    \allowdisplaybreaks
\usepackage{amssymb}
\usepackage{bm}
\usepackage{physics}
\usepackage[mathscr]{euscript}
\usepackage{multirow}
\usepackage[inline]{enumitem}
\usepackage[x11names]{xcolor}
\usepackage{hyperref}
    \hypersetup{
        linkbordercolor=Goldenrod1,
        citebordercolor=LightSkyBlue1,
        urlbordercolor=OliveDrab1,
    }
\usepackage{amsthm}
    \theoremstyle{plain} 
        \newtheorem{theorem}{Theorem}
        
        \newtheorem{lemma}{Lemma}
        \newtheorem{proposition}{Proposition}
    \theoremstyle{definition}
        \newtheorem{assumption}{Assumption}
        \newtheorem{definition}{Definition}
        \newtheorem{example}{Example}
    \theoremstyle{remark}
        \newtheorem{remark}{Remark}
\usepackage{cleveref}
    
    \crefname{assumption}{assumption}{assumptions}
\usepackage[margin={1in,1in}]{geometry}
\usepackage{booktabs}
\usepackage{microtype}
\usepackage{setspace}
\usepackage{algorithm}
\usepackage{algpseudocode}
\usepackage[labelfont=bf]{caption}
\usepackage{natbib}
    \bibpunct[, ]{(}{)}{,}{a}{}{,}%
\usepackage{gentium}
\usepackage{marvosym} 

\makeatletter
\def\blfootnote{\gdef\@thefnmark{}\@footnotetext}
\makeatother

\input{macros}

\DeclareMathOperator{\KG}{KG}
\DeclareMathOperator{\IKG}{IKG}
\def\acc{\mathsf{acc}}
\def\RUNTITLE{Knowledge Gradient for Selection with Covariates: Consistency and Computation} 
\def\RUNAUTHOR{Ding et al.} 

\usepackage[headsepline]{scrlayer-scrpage}
    \ohead{\pagemark}
    \ihead{\footnotesize \emph{\textbf{\RUNAUTHOR}}: \RUNTITLE}
    \ofoot[]{}

\title{Knowledge Gradient for Selection with Covariates: Consistency and Computation}

\author{
    {Liang Ding\thanks{
        Department of Industrial Systems Engineering, Texas A\&M University, College Station, TX, USA
        }
    }
    \and
    {L. Jeff Hong\thanks{
        School of Management and School of Data Science, Fudan University, Shanghai, China
        }
    }
    \and
    {Haihui Shen\thanks{
        Sino-US Global Logistics Institute, Antai College of Economics and Management, Shanghai Jiao Tong University, Shanghai, China
        }
    \hspace{1.8pt}$^\text{\Letter}$
    }
    \and
    {Xiaowei Zhang\thanks{
        Faculty of Business and Economics, University of Hong Kong, Pok Fu Lam, Hong Kong SAR
        }
    }
}

\date{}

\begin{document}

\pdfbookmark[1]{Title}{title}

\maketitle

\vspace{-10pt}
\begin{abstract}
    \noindent
    Knowledge gradient is a design principle for developing Bayesian sequential sampling policies to solve optimization problems.
    In this paper we consider the ranking and selection problem in the presence of covariates, where the best alternative is not universal but depends on the covariates.
    In this context, we prove that under minimal assumptions, the sampling policy based on knowledge gradient is consistent, in the sense that following the policy the best alternative as a function of the covariates will be identified almost surely as the number of samples grows.
    We also propose a stochastic gradient ascent algorithm for computing the sampling policy and demonstrate its performance via numerical experiments.

\end{abstract}
\vspace{10pt}

\blfootnote{\hspace{-5pt}$^\text{\Letter}$\hspace{-0.5pt}Corresponding Author, \href{mailto:shenhaihui@sjtu.edu.cn}{shenhaihui@sjtu.edu.cn}}

\section{Introduction}

We consider the ranking and selection (R\&S) problem in the presence of covariates.
A decision maker is presented with a finite collection of alternatives. The performance of each alternative is unknown and depends on the covariates.
Suppose that the decision maker has access to noisy samples of each alternative for any chosen value of the covariates, but the samples are expensive to acquire.
Given a finite sampling budget, the goal is to develop an efficient sampling policy indicating locations as to which alternative and what value of the covariates to sample from, so that upon termination of the sampling, the decision maker can identify a decision rule that accurately specifies the best alternative as a function of the covariates.

The problem of R\&S with covariates emerges naturally as the popularization of data and decision analytics in recent years.
In clinical and medical research, for many diseases the effect of a treatment may be substantially different across patients, depending on their biometric characteristics (i.e., the covariates), including
age, weight, lifestyle habits such as smoking and alcohol use, etc.~\citep{kim2011short}. A treatment regime that works for a majority of patients might not work for the others.
Samples needed for estimating treatment effects may be collected from clinical trials or computer simulation.
For example, in \cite{hur2004} and \cite{choi2014}, a simulation model is developed to simulate the effect of several treatment regimens for Barrett's esophagus, a precursor to esophageal cancer, for patients with different biometric characteristics.
Personalized medicine can then be developed to determine the best treatment regime that is customized to the particular characteristics of each individual patient.
Similar customized decision-making can be found in online advertising~\citep{arora2008}, where advertisements are displayed depending on consumers' web browsing history or buying behavior to increase the revenue of the advertising platform as well as to improve consumers' shopping experience.

Being a classic problem in the area of stochastic simulation, R\&S has a vast literature. We refer to~\cite{kim2006} and~\cite{ChenChickLeePujowidianto15} for reviews on the subject with emphasis on frequentist and Bayesian approaches, respectively. Most of the prior work, however, does not consider the presence of the covariates, and thus the best alternative to select is universal rather than varies as a function of the covariates. There are several exceptions, including~\cite{HuLudkovski17}, ~\cite{PearceBranke17}, and \cite{ShenHongZhang17}. Among them \cite{ShenHongZhang17} take a frequentist approach to solve R\&S with covariates, whereas the other two a Bayesian approach. The present paper adopts a Bayesian perspective as well.

This paper considers a sampling policy based on knowledge gradient (KG) for R\&S with covariates.
KG, introduced in~\cite{FrazierPowellDayanik08}, is a design principle that has been widely used for developing Bayesian sequential sampling policies to solve a variety of optimization problems, including R\&S, in which evaluation of the objective function is noisy and expensive. In its basic form, KG begins with assigning a multivariate normal prior on the unknown constant performance of all alternatives. In each iteration, it chooses the sampling location by maximizing the increment in the expected value of the information that would be gained by taking a sample from the location. Then, the posterior is updated upon observing the noisy sample from the chosen location. The sampling efficiency of KG-type policies is often competitive with or outperforms other sampling policies; see~\cite{FrazierPowellDayanik09},~\cite{ScottFrazierPowell11},~\cite{Ryzhov16}, and~\cite{PearceBranke18} among others.

A KG-based sampling policy for R\&S with covariates is also proposed in ~\cite{PearceBranke17}. The main difference here is that our treatment is more general. First, we allow the sampling noise to be heteroscedastic, whereas it is assumed to be constant for different locations of the same alternative in their work.
Heteroscedasticity is of particular significance for simulation applications such as queueing systems.
Second, we take into account possible variations in sampling cost at different locations, whereas the sampling cost is simply treated as constant everywhere in~\cite{PearceBranke17}. Hence, our policy, which we refer to as integrated knowledge gradient (IKG), attempts in each iteration to maximize a ``cost-adjusted'' increment in the expected value of information.
These generalizations are straightforward when the variance of the sampling noise and the sampling cost are assumed to be known.
We also briefly discuss and show how to deal with the case where they are unknown.

The first main contribution of this paper is to provide a theoretical analysis of the asymptotic behavior of the IKG policy, whereas~\cite{PearceBranke17}  conducted only numerical investigation. In particular, we prove that IKG is consistent in the sense that for any value of the covariates, the selected alternative upon termination of the sampling will converge to the true best almost surely as the sampling budget grows to infinity.
Moreover, we consider a practical variant---termed quasi-IKG---which does not require the intermediate optimization problem in each iteration of IKG to be solved exactly,
and prove its consistency under mild conditions.

Consistency of KG-type policies has been established in various settings, mostly for problems where the number of feasible solutions is finite, including R\&S~\citep{FrazierPowellDayanik08,FrazierPowellDayanik09,FrazierPowell11,MesPowellFrazier11}, and discrete optimization via simulation~\citep{XieFrazierChick16}. KG is also used for Bayesian optimization of continuous functions in~\cite{WuFrazier16},~\cite{PoloczekWangFrazier17}, and~\cite{WuPoloczekWilsonFrazier17}. However, in these papers the continuous domain is discretized first, which effectively reduces the problem to one with finite feasible solutions, in order to facilitate their asymptotic analysis.
The finiteness of the domain is critical in the aforementioned papers, because the asymptotic analysis there boils down to proving that each feasible solution can be sampled infinitely often. This, by the law of large numbers, implies that the variance of the objective value estimate of each solution will converge to zero. Thus, the optimal solution will be identified ultimately since the uncertainty about the performances of the solutions will be removed completely in the end.

By contrast, proving consistency of KG-type policies for continuous solution domains demands a fundamentally different approach, since most solutions in a continuous domain would hardly be sampled even once after all.
Among the several related papers, \cite{ScottFrazierPowell11} studies a KG-type policy for Bayesian optimization of continuous functions.
Assigning a Gaussian process prior on the objective function, they established the consistency of the KG-type policy basically by leveraging the continuity of the covariance function of the Gaussian process, which intuitively suggests that if the variance at one location is small, then the variance in its neighborhood ought to be small too.
\cite{toscano2018bayesian} prove the consistency of a KG-type policy on a more general problem that can reduce to the problem in \cite{ScottFrazierPowell11}, for both discrete and continuous domains.

We cast R\&S with covariates to a problem of ranking a finite number of Gaussian processes, thereby having both discrete and continuous elements structurally. As a result, we establish the consistency of the proposed IKG policy by proving the following two facts -- (i) each Gaussian process is sampled infinitely often, and (ii) the infinitely many samples assigned to a given Gaussian process drives its posterior variance at any location to zero, thanks to the assumed continuity of its covariance function. The theoretical analysis in this paper is partly built on the ideas developed for discrete and continuous problems, respectively, in~\cite{FrazierPowellDayanik08} and~\cite{ScottFrazierPowell11} in a federated manner.

Although our proofs share similar structures to those in \cite{ScottFrazierPowell11}, our assumptions are substantially simpler and minimal.
By contrast, for the proof in~\cite{ScottFrazierPowell11} to be valid, technical conditions are imposed to regulate the asymptotic behavior of the posterior mean function and the posterior covariance function of the underlying Gaussian process. Nevertheless, the two conditions are difficult to verify. We do not impose such conditions. We achieve the substantial simplification of the assumptions by leveraging the reproducing kernel Hilbert space (RKHS) theory. The theory has been used widely in machine learning~\citep{SteinwartChristmann08}. But its use in the analysis of KG-type policies is less common.
We develop several technical results based on RKHS theory to facilitate analysis of the asymptotic behavior of the posterior covariance function.\footnote{\citet{bect2019} adopt a  supermartingale approach to study the asymptotic behavior of a general class of sequential sampling algorithms.
Their analysis has a broader scope of applicability but it is technically more involved.}

The second
main contribution of this paper is that we develop an algorithm to solve a stochastic  optimization problem that determines the sampling decision of the IKG policy in its each iteration. In~\cite{PearceBranke17}, this optimization problem is addressed by the sample average approximation method with a derivative-free optimization solver. Instead, we propose a stochastic gradient ascent (SGA) algorithm, taking advantage of the fact that a gradient estimator can be derived analytically for many popular covariance functions. Numerical experiments demonstrate the finite-sample performance of the IKG policy in conjunction with the SGA algorithm.

We conclude the introduction by reviewing briefly the most pertinent literature. A closely related problem is multi-armed bandit (MAB); see~\cite{bubeck2012} for a comprehensive review on the subject. The significance of covariates, thereby contextual MAB (or MAB with covariates), has also drawn substantial attention in recent years; see~\cite{rusmevichientong2010},~\cite{yang2002},~\cite{KrauseOng11}, and~\cite{perchet2013} among others. There are two critical differences between contextual MAB and R\&S with covariates. First, the former generally assumes that the covariates arrive exogenously in a sequential manner, and the decision-maker can choose  at which arm (or alternative) to sample but not the value of covariates. By contrast, the latter assumes that the decision-maker is capable of choosing both the alternative and the covariates when specifying sampling locations. A second difference is MAB focuses on minimizing the regret which is caused by choosing inferior alternatives and accumulated during the sampling process, whereas R\&S focuses on identifying the best alternative eventually and the regret is not the primary concern.

The rest of the paper is organized as follows.
In  \Cref{sec:formulation} we follow a nonparametric Bayesian approach to formulate the problem of R\&S with covariates, introduce the IKG policy, and present the main result.
In \Cref{sec:consistency} we prove the consistency of our sampling policy in the sense that the estimated best alternative as a function of the covariates converges to the truth with probability one as the number of samples grows to infinity.
We then propose to use SGA for computing our sampling policy in \Cref{sec:computation}, and demonstrate its performance via numerical experiments in \Cref{sec:numerical}.
We conclude in \Cref{sec:conclusion} and collect detailed proof and additional technical results and numerical experiments in the Appendix.

\section{Problem Formulation}\label{sec:formulation}

Suppose that a decision maker is presented with \(M\) competing alternatives.
For each \(i=1,\ldots,M\), the performance of alternative \(i\) depends on a vector of \emph{covariates} \(\BFx=(x_1,\ldots,x_d)^\intercal\) and is denoted by \(\theta_i=\theta_i(\BFx)\) for \(\BFx\in\calX\subset \Real^d\).
The performances are unknown and can only be learned via sampling.
In particular, for any \(i\) and \(\BFx\), one can acquire possibly multiple noisy samples of \(\theta_i(\BFx)\).
The decision maker aims to select the ``best'' alternative for a given value of \(\BFx\), i.e., identify  \(\argmax_i\theta_i(\BFx)\).
However, since the sampling is usually expensive in time and/or money, instead of estimating the performances \(\{\theta_i(\BFx):i=1,\ldots,M\}\) every time a new value of \(\BFx\) is observed and then ranking them, it is preferable to learn \emph{offline} the decision rule
\begin{equation}\label{eq:decision_rule}
    i^*(\BFx) \in \argmax_{1\leq i\leq M} \theta_i(\BFx), \quad \BFx\in\calX,
\end{equation}
as a function of \(\BFx\), through a carefully designed sampling process.
Equipped with such a decision rule, the decision maker can select the best alternative upon observing the covariates in a timely fashion.
In addition, the decision maker may have some knowledge with regard to the covariates. For example, certain values of the covariates may be more important or appear more frequently than others.  Suppose that this kind of knowledge is expressed by a probability density function \(\gamma(\BFx)\) on \(\calX\).

During the offline learning period,
we need to make a sequence of sampling decisions \(\{(a^n,\BFv^n):n=0,1,\ldots\}\), where \((a^n,\BFv^n)\) means that the \((n+1)\)-th sample, denoted by \(y^{n+1}\), is taken from alternative \(a^n\) with covariates value \(\BFv^n\) (refer it as location \(\BFv^n\) for simplicity).
We assume that given \(\theta_{a^n}(\BFv^n)\), \(y^{n+1}\) is an unbiased sample having a normal distribution, i.e.,
\begin{equation*}\label{eq:dist_of_sample}
    y^{n+1}\,|\,\theta_{a^n}(\BFv^n) \sim \calN(\theta_{a^n}(\BFv^n), \lambda_{a^n}(\BFv^n)),
\end{equation*}
where $y^n\,|\,\theta_i(\BFx)$ is independent of $y^{n'}\,|\,\theta_{i'}(\BFx')$ for $(i,\BFx,n) \neq (i',\BFx',n')$.
Here, \(\lambda_i(\BFx)\) is the variance of a sample of \(\theta_i(\BFx)\) given \(\theta_i(\BFx)\) and is assumed to be \emph{known}.
Moreover, suppose that the cost of taking a sample from alternative \(i\) at location \(\BFx\) is \(c_i(\BFx) >0\), which is also assumed to be \emph{known}.
Suppose that the total sampling budget for offline learning is \(B>0\), and the sampling process is terminated when the budget is exhausted.
Mathematically, we will stop with the \(N(B)\)-th sample, where
\begin{equation}\label{eq:budget}
    N(B) \coloneqq \min \qty{N: \sum_{n=0}^N c_{a^n}(\BFv^n) > B }.
\end{equation}
Consequently, the sampling decisions are \(\{(a^n,\BFv^n):n=0,\ldots, N(B)-1\}\) and the samples taken during the process are \(\{y^{n+1}:n=0,\ldots, N(B)-1\}\).
Notice that \(N(B)=B\) if \(c_i(\BFx)\equiv1\) for \(i=1,\ldots,M\), in which case the sampling budget is reduced to the number of samples.

\begin{remark}
The assumption of known \(\lambda_i(\BFx)\) is critical to the theoretical analysis in this paper.
As we will see shortly, with known \(\lambda_i(\BFx)\), if we impose a Gaussian process as prior for $\theta_i$, then its posterior will still be a Gaussian process, which makes the asymptotic analysis tractable.
It would not be the case if \(\lambda_i(\BFx)\) also needs to be estimated.
In practice, \(\lambda_i(\BFx)\) is usually unknown and it is a common issue in the experiment design.
We suggest to follow the approach in \cite{AnkenmanNelsonStaum10}, which fits the surfaces of \(\lambda_i(\BFx)\) by running multiple simulations at certain design points and using the sample variances.
See more details in the numerical experiments in the Appendix.
The unknown sampling cost \(c_i(\BFx)\) in practice can be handled similarly.
\end{remark}

We follow a nonparametric Bayesian approach to model the unknown functions \(\{\theta_1,\ldots,\theta_M\}\) as well as to design the sampling policy.
We treat \(\theta_i\)'s as random functions and impose a prior on them under which they are mutually independent, although this assumption may be relaxed.
Suppose that \(\BFx\) takes continuous values and that under the prior, \(\theta_i\) is a Gaussian process with mean function \(\mu_i^0(\BFx)\coloneqq\E[\theta_i(\BFx)]\) and covariance function \(k_i^0(\BFx,\BFx')\coloneqq\Cov[\theta_i(\BFx), \theta_i(\BFx')]\) that satisfies the following assumption.

\begin{assumption}\label{assump:stationary}
    For each \(i=1,\ldots,M\), there exists a constant \(\tau_i> 0\) and a positive continuous function \(\rho_i:\Real^d\mapsto\Real_{+}\) such that  \(k_i^0(\BFx,\BFx')=\tau_i^2\rho_i(\BFx-\BFx')\).
    Moreover,
    \begin{enumerate}[label=(\roman*)]
        \item
              \(\rho_i(\abs{\BFdelta})=\rho_i(\BFdelta)\), where \(\abs{\cdot}\) means taking the absolute value component-wise;
        \item
              \(\rho_i(\BFdelta)\) is  decreasing in \(\BFdelta\) component-wise for \(\BFdelta\geq \BFzero\);
        \item
              \(\rho_i(\BFzero)= 1\), \(\rho_i(\BFdelta)\to 0\) as \(\norm{\BFdelta}\to\infty\), where \(\norm{\cdot}\) denotes the Euclidean norm;
        \item
              there exist some $0<C_i<\infty$ and $\varepsilon_i, u_i>0$ such that
              $$1-\rho_i(\BFdelta)\leq \frac{C_i} {|\log(\norm{\BFdelta})|^{1+\varepsilon_i}},$$
              for all $\BFdelta$ such that $\norm{\BFdelta} < u_i$.
    \end{enumerate}
\end{assumption}

\begin{remark}
\Cref{assump:stationary} stipulates that \(k_i^0\) is
\emph{stationary}, i.e., it depends on \(\BFx\) and \(\BFx'\) only through the difference \(\BFx-\BFx'\).
In addition, \(\tau_i^2\) can be interpreted as the prior variance of \(\theta_i(\BFx)\) for all \(\BFx\), and \(\rho_i(\BFx-\BFx')\) as the prior correlation between \(\theta_i(\BFx)\) and \(\theta_i(\BFx')\) which increases to 1 as \(\norm{\BFx-\BFx'}\) decreases to 0.
The condition in part (iv) of \Cref{assump:stationary} is weak.
In conjunction with the continuity assumption of $\rho_i$, it implies that the sample paths of the Gaussian process $\theta_i$ are continuous almost surely if the mean function $\mu_i^0(\BFx)$ is continuous; see, e.g., \citet[Theorem 1.4.1]{AdlerTaylor07}.
The sample path continuity will be used to establish the uniform convergence of the posterior mean functions.
\end{remark}

A variety of covariance functions satisfy \Cref{assump:stationary}.
Notable examples include the squared exponential (SE) covariance function
\[k_{\text{SE}}(\BFx,\BFx')=\tau^2 \exp\qty(- r^2(\BFx-\BFx')),\]
where \(r(\BFdelta)=\sqrt{\sum_{j=1}^d\alpha_j \delta_j^2}\) and \(\alpha_j\)'s are positive parameters, and the Mat\'{e}rn covariance function
\[k_{\text{Mat\'{e}rn}}(\BFx,\BFx')=\tau^2 \frac{2^{1-\nu}}{\Gamma(\nu)}\qty(\sqrt{2\nu}r(\BFx-\BFx'))^\nu K_\nu\qty(\sqrt{2\nu}r(\BFx-\BFx')),\]
where \(\nu\) is a positive parameter that is typically taken as half-integer (i.e., \(\nu=p+1/2\) for some nonnegative integer \(p\)), \(\Gamma\) is the gamma function, and \(K_\nu\) is the modified Bessel function of the second kind.
The covariance function reflects one's prior belief about the unknown functions. We refer to~\citet[Chapter 4]{RasmussenWilliams06} for more types of covariance functions.

\subsection{Bayesian Updating Equations}

For each \(n=1,2,\ldots\), let \(\mathscr{F}^n\) denote the \(\sigma\)-algebra generated by \((a^0, \BFv^0),y^1,\ldots,(a^{n-1},\BFv^{n-1}),y^n\), the sampling decisions and the samples collected up to time \(n\).
Suppose that \((a^n,\BFv^n)\in\mathscr{F}^n\), that is, \((a^n,\BFv^n)\) depends only on the information available at time \(n\).
In addition, we use the notation \(\E^n[\cdot]\coloneqq \E[\cdot|\mathscr{F}^n]\), and define \(\Var^n[\cdot]\) and \(\Cov^n[\cdot]\) likewise.

Given the setup of our model, it is easy to derive that \(\{\theta_1,\ldots,\theta_M\}\) are independent Gaussian processes under the posterior distribution conditioned on \(\mathscr F^n\), \(n=1,\ldots,N(B)\).
In particular, under the prior mutual independence, taking samples from one unknown function does not provide information on another.
Let \(\BFV_i^n\coloneqq\{\BFv^\ell: a^\ell = i, \ell=0,\ldots,n-1\}\) denote  the set of the locations of the samples taken from \(\theta_i\) up to time \(n\) and define \(\BFy_i^n \coloneqq\{y^{\ell+1}: a^\ell = i, \ell=0,\ldots,n-1\}\) likewise.
With slight abuse of notation, when necessary, we will also treat \(\BFV_i^n\) as a matrix wherein the columns are corresponding to the points in the set and arranged in the order of appearance, and \(\BFy_i^n\) as a column vector with elements also arranged in the order of appearance.
Then, the posterior mean and covariance functions of \(\theta_i\) are given by
\begin{align}
    \mu_i^n(\BFx)     & \coloneqq \E^n[\theta_i(\BFx)]  = \mu_i^0(\BFx) +  k_i^0(\BFx,\BFV_i^n)[k_i^0(\BFV_i^n,\BFV_i^n)+\lambda_i(\BFV_i^n)]^{-1}[\BFy_i^n-\mu_i^0(\BFV_i^n)], \label{eq:mean_n}             \\[0.5ex]
    k_i^n(\BFx,\BFx') & \coloneqq \Cov^n[\theta_i(\BFx), \theta_i(\BFx')] = k_i^0(\BFx,\BFx') - k_i^0(\BFx,\BFV_i^n)[k_i^0(\BFV_i^n,\BFV_i^n)+\lambda_i(\BFV_i^n)]^{-1}k_i^0(\BFV_i^n,\BFx'),\label{eq:cov_n}
\end{align}
where for two sets \(\BFV\) and \(\BFV'\), \(k_i^0(\BFV,\BFV')=[k_i^0(\BFx,\BFx')]_{\BFx\in\BFV,\BFx'\in\BFV'}\) is a matrix of size \(\abs{\BFV}\times \abs{\BFV'}\), \(\lambda_i(\BFV)=\mathrm{diag}\{\lambda_i(\BFx):\BFx\in\BFV\}\) is a diagonal matrix of size \(\abs{\BFV}\times \abs{\BFV}\), and \(\mu_i^0(\BFV)=(\mu_i^0(\BFx):\BFx\in\BFV)\) is a column vector of size \(\abs{\BFV}\times 1\).
Here $|\cdot|$ denotes the cardinality of a set.
We refer to, for example,~\citet[Section 3.2]{ScottFrazierPowell11} for details.
Further, the following updating equation can be derived
\begin{align}
    \mu_i^{n+1}(\BFx)     & = \mu_i^n(\BFx) + \sigma_i^n(\BFx,\BFv^n) Z^{n+1},\label{eq:mean_update}                    \\[0.5ex]
    k_i^{n+1}(\BFx,\BFx') & = k_i^n(\BFx,\BFx') - \sigma_i^n(\BFx,\BFv^n)\sigma_i^n(\BFx',\BFv^n),\label{eq:cov_update}
\end{align}
where \(Z^{n+1}\) is a standard normal random variable independent to everything else, and
\begin{equation}\label{eq:sigma}
    \sigma_i^n(\BFx,\BFv^n) \coloneqq
    \begin{cases}
        \tilde\sigma_i^n(\BFx,\BFv^n), & \text{if } i = a^n,    \\
        0,                             & \text{if }  i\neq a^n,
    \end{cases}
    \qq{and}
    \tilde\sigma_i^n(\BFx, \BFv) \coloneqq \frac{k_i^n(\BFx,\BFv)}{\sqrt{k_i^n(\BFv,\BFv)+\lambda_i(\BFv)}}.
\end{equation}
In particular, conditioned on \(\mathscr{F}^n\) and prior to taking a sample at \((a^n,\BFv^n)\), the predictive distribution of \(\mu_i^{n+1}(\BFx)\) is normal with mean \(\mu_i^n(\BFx)\) and standard deviation \(\sigma_i^n(\BFx,\BFv^n)\).
Moreover, notice that
\begin{equation}\label{eq:decreasing_var}
    \Var^{n+1}[\theta_i(\BFx)] = k_i^{n+1}(\BFx,\BFx) = k_i^n(\BFx,\BFx) - [\sigma_i^n(\BFx,\BFv^n)]^2\leq \Var^n[\theta_i(\BFx)].
\end{equation}
(Note that \crefrange{eq:mean_update}{eq:decreasing_var} are still valid even if \(k_i^n(\BFv^n,\BFv^n)=0\), and/or \(\lambda_i(\BFv^n)=0\).)
Hence, \(\Var^n[\theta_i(\BFx)]\) is non-increasing in \(n\).
This basically suggests that the uncertainty about each unknown function under the posterior decreases as more samples from it are collected.
It is thus both desirable and practically meaningful that such uncertainty would be completely eliminated if the sampling budget is unlimited, in which case one would be able to identify the decision rule \cref{eq:decision_rule} perfectly.
In particular, we define consistency of a sampling policy as follows.

\begin{definition}
    A sampling policy is said to be \emph{consistent} if it ensures that
    \begin{equation}\label{eq:consistent}
        \lim_{B\to\infty}\argmax_{1\leq i\leq M}\mu^{N(B)}_i(\BFx)= \argmax_{1\leq i\leq M}\theta_i(\BFx),
    \end{equation}
    almost surely (a.s.) for all \(\BFx\in\calX\).
\end{definition}

\begin{remark}
    Under the assumption that \(\{\theta_1,\ldots,\theta_M\}\) are independent under the prior, collecting samples from \(\theta_i\) does not provide information about \(\theta_j\) if \(i\neq j\).
    Therefore, a consistent policy under the independence assumption ought to ensure that the number of samples taken from \emph{each} \(\theta_i\) grows  without bounds.
\end{remark}

\subsection{Knowledge Gradient Policy}

We first assume temporarily that \(\BFx\) is given and fixed,
and that \(c_i(\BFx)= 1\) for \(i=1,\ldots,M\).
Then, solving \(\max_i \theta_i(\BFx)\) is a selection of the best problem having finite alternatives, and each sampling decision is reduced to choosing an alternative \(i\) to take a sample of \(\theta_i(\BFx)\).
The knowledge gradient (KG) policy introduced in~\cite{FrazierPowellDayanik08} is designed exactly to solve such a problem assuming an independent normal prior.
Specifically, the knowledge gradient at \(i\) is defined there as the increment in the expected value of the information about the maximum at \(\BFx\) gained by taking a sample at \(i\), that is,
\begin{equation}\label{eq:KG}
    \KG^n(i;\BFx)\coloneqq \E\qty[\max_{1\leq a\leq M} \mu^{n+1}_a(\BFx)\,\Big|\,\mathscr{F}^n, a^n=i] - \max_{1\leq a\leq M} \mu^n_a(\BFx).
\end{equation}
Then, each time the alternative \(i\) that has the largest value of \(\KG(i;\BFx)\) is selected to generate a sample of \(\theta_i(\BFx)\).

Let us now return to our context where (1) the covariates are present, (2) each sampling decision consists of both \(i\) and \(\BFx\), and (3) each sampling decision may induce a different sampling cost.
Since a sample of \(\theta_i(\BFx)\) would alter the posterior belief about \(\theta_i(\BFx')\), we generalize \cref{eq:KG} and define
\begin{equation}\label{eq:new_KG}
    \KG^n(i,\BFx;\BFv)\coloneqq \frac{1}{c_i(\BFx)} \qty{ \E\qty[\max_{1\leq a\leq M} \mu^{n+1}_a(\BFv)\,\Big|\,\mathscr{F}^n, a^n=i, \BFv^n=\BFx] - \max_{1\leq a\leq M} \mu^n_a(\BFv)},
\end{equation}
which can be interpreted as the increment in the expected value of the information about the maximum at \(\BFv\) gained per unit of sampling cost
by taking a sample at \((i,\BFx)\).
Then, we consider the following \emph{integrated KG} (IKG)
\begin{equation}\label{eq:int-KG}
    \IKG^n(i,\BFx)\coloneqq \frac{1}{c_i(\BFx)} \int_{\calX}\qty{\E\qty[\max_{1\leq a\leq M} \mu^{n+1}_a(\BFv)\,\Big|\,\mathscr{F}^n, a^n=i, \BFv^n=\BFx] - \max_{1\leq a\leq M} \mu^n_a(\BFv)} \gamma(\BFv) \dd{\BFv},
\end{equation}
and define the IKG sampling policy as
\begin{equation}\label{eq:max-int-KG}
    (a^n,\BFv^n) \in \argmax_{1\leq i\leq M, \BFx\in\calX}\IKG^n(i,\BFx).
\end{equation}
The integrand of \cref{eq:int-KG} can be calculated analytically, as shown in \Cref{lemma:calculating_IKG}, whose proof is deferred to the Appendix.

\begin{lemma}\label{lemma:calculating_IKG}
    For all \(i=1,\ldots,M\) and \(\BFx\in\calX\),
    \begin{equation}\label{eq:int-KG2}
        \IKG^n(i,\BFx)= \frac{1}{c_i(\BFx)} \int_{\calX}\left[\abs{\tilde\sigma_i^n(\BFv, \BFx)}\phi\qty(\abs{\frac{\Delta_i^n(\BFv)}{\tilde\sigma_i^n(\BFv, \BFx)}})
            -\abs{\Delta_i^n(\BFv)}\Phi\qty(-\abs{\frac{\Delta_i^n(\BFv)}{\tilde\sigma_i^n(\BFv, \BFx)}})
            \right] \gamma(\BFv) \dd\BFv,
    \end{equation}
    where \(\Delta_i^n(\BFv)\coloneqq \mu_i^n(\BFv)-\max_{a\neq i}\mu_a^n(\BFv)\), \(\Phi\) is the standard normal distribution function, and \(\phi\) is its density function.
\end{lemma}

We solve \cref{eq:max-int-KG} by first solving \(\max_{\BFx}\IKG^n(i,\BFx)\) for all \(i\) and then enumerating the results.
The computational challenge in the former lies in the numerical integration in \cref{eq:int-KG2}.
Notice that \(\max_{\BFx}\IKG^n(i,\BFx)\)  is in fact a stochastic optimization problem if we view the integration in \cref{eq:int-KG2} as an expectation with respect to the probability density \(\gamma(\BFx)\) on \(\calX\).
One might apply  the sample average approximation method to solve \(\max_{\BFx}\IKG^n(i,\BFx)\), but it would be computationally prohibitive if \(\calX\) is high-dimensional.
Instead, we show in \Cref{sec:computation} that the gradient of the integrand in \cref{eq:int-KG2} with respect to \(\BFx\) can be calculated explicitly, which is an unbiased estimator of \(\nabla_{\BFx}\IKG^n(i,\BFx)\) under regularity conditions, thereby leading to a stochastic gradient ascent method~\citep{KushnerGeorge03}.

We now present our main theoretical result --- the IKG policy is consistent under simple assumptions.
The proof will be sketched in \Cref{sec:consistency} and all details are collected in the Appendix.

\begin{assumption}\label{assump:compact}
    The design space \(\calX\) is a compact set in \(\Real^d\) with nonempty interior.
\end{assumption}

\begin{assumption}\label{assump:var}
    For each \(i=1,\ldots,M\), \(\mu_i^0(\cdot)\), \(\lambda_i(\cdot)>0\) and \(c_i(\cdot)>0\) are all continuous on \(\calX\), and  \(\gamma(\cdot) >0\) on \(\calX\).
\end{assumption}

Under \Cref{assump:stationary,assump:compact,assump:var}, the IKG policy \eqref{eq:max-int-KG} is well defined.
This can be seen by noting that the maximum of \(\IKG^n(i,\BFx)\) over \(\BFx\in\calX\) is attainable since \(\IKG^n(i,\BFx)\) is continuous in \(\BFx\) by \Cref{assump:stationary,assump:var} together with \Cref{lemma:calculating_IKG}, and \(\calX\) is compact by \Cref{assump:compact}.
Moreover, the IKG policy \eqref{eq:max-int-KG} is consistent as formally stated in the following \Cref{theo:consistency-cost}.

\begin{theorem}\label{theo:consistency-cost}
    If \Cref{assump:stationary,assump:compact,assump:var} hold, then the IKG policy \eqref{eq:max-int-KG} is  consistent, that is, under the IKG policy,
    \begin{enumerate}[label=(\roman*)]
        \item
              \(k_i^{N(B)}(\BFx,\BFx)\to 0\) a.s. as \(B\to\infty\) for all \(\BFx\in\calX\) and \(i=1,\ldots,M\);
        \item
              \(\mu^{N(B)}_i(\BFx)\to \theta_i(\BFx)\) a.s. as \(B\to\infty\)  for all \(\BFx\in\calX\) and \(i=1,\ldots,M\);
        \item
              \(\argmax_{1\leq i\leq M}\mu^{N(B)}_i(\BFx)\to \argmax_{1\leq i\leq M}\theta_i(\BFx)\) a.s.  as \(B\to\infty\) for all \(\BFx\in\calX\).
    \end{enumerate}
\end{theorem}

We conclude this section by highlighting the differences between our assumptions and those in~\cite{ScottFrazierPowell11}, in which the consistency of a KG-type policy driven by a Gaussian process is proved. First and foremost, \emph{they impose conditions on both the posterior mean function and the posterior covariance function to regulate their large-sample asymptotic behavior}.
Specifically, they assume that  uniformly for all \(n\) and \(\BFx,\BFv\in\calX\) with \(\BFx\neq \BFv\), (1) \(\abs{\mu^n(\BFx)-\mu^n(\BFv)}\) is bounded a.s., and (2) \(\abs{\Corr^n[\theta(\BFx),\theta(\BFv)]}\) is bounded above away from one, where \(\Corr^n\) means the posterior correlation.\footnote{The subscript \(i\) is ignored because there is only one Gaussian process involved in~\cite{ScottFrazierPowell11}.}
\emph{The two assumptions are critical for their analysis but nontrivial to verify in practice}.

By contrast, we do not make such assumptions.
Condition (1) is not necessary in our analysis because the ``increment in the expected value of the information'' is defined as \cref{eq:int-KG} in this paper, whereas in a different form without integration in~\cite{ScottFrazierPowell11}. There is no need for us to impose Condition (2) in order to regulate  the asymptotic behavior of the posterior covariance function, because instead we achieve the same goal by utilizing reproducing kernel Hilbert space (RKHS) theory.

Second, in~\cite{ScottFrazierPowell11} the prior covariance function of the underlying Gaussian process is of SE type. We relax it to \Cref{assump:stationary}, which allows a great variety of covariance functions.
We also take into account possibly varying sampling costs at different locations.

\section{Consistency}\label{sec:consistency}

It is straightforward to show that \(N(B) \to \infty\) if and only if \(B \to \infty\), since \(c_i(\cdot)\) is bounded both above and below away from zero on \(\calX\) for each \(i=1,\ldots,M\) under \Cref{assump:compact,assump:var}. Thus, \Cref{theo:consistency-cost} is equivalent to \Cref{theo:consistency} as follows.

\begin{theorem}\label{theo:consistency}
    If \Cref{assump:stationary,assump:var,assump:compact} hold, then under the IKG policy,
    \begin{enumerate}[label=(\roman*)]
        \item
              \(k_i^n(\BFx,\BFx)\to 0\) a.s. as \(n\to\infty\) for all \(\BFx\in\calX\) and \(i=1,\ldots,M\);
        \item
              \(\mu^n_i(\BFx)\to \theta_i(\BFx)\) a.s. as \(n\to\infty\)  for all \(\BFx\in\calX\) and \(i=1,\ldots,M\);
        \item
              \(\argmax_{1\leq i\leq M}\mu^n_i(\BFx)\to \argmax_{1\leq i\leq M}\theta_i(\BFx)\) a.s.  as \(n\to\infty\) for all \(\BFx\in\calX\).
    \end{enumerate}
\end{theorem}

The bulk of the proof of consistency of the IKG policy lies in part (i) of \Cref{theo:consistency}, i.e., to show that \(\lim_{n\to\infty}\Var^n[\theta_i(\BFx)]= 0\) a.s. for all \(\BFx\in\calX\) and \(i=1,\ldots,M\).
It consists of two steps, which are summarized into the later \Cref{prop:step1,prop:step2}.
However, both \Cref{prop:step1,prop:step2} critically relies on the asymptotic behavior of the posterior covariance function, which is characterized in the following \Cref{prop:unif_conv}.

\begin{proposition}\label{prop:unif_conv}
    Fix \(i=1,\ldots,M\).
    If \(k^0_i\) is stationary, then for any \(\BFx\in\calX\), \(k_i^n(\BFx,\BFx')\) converges to a limit, denoted by \(k_i^\infty(\BFx,\BFx')\), uniformly in \(\BFx'\in\calX\) as \(n\to\infty\).
\end{proposition}

\Cref{prop:unif_conv} shows that \emph{irrespective} of the allocation of the design points \(\{\BFv^\ell:\ell=0,\ldots,n-1\}\), \(k_i^n(\BFx,\cdot)\) converges \emph{uniformly} as \(n\to\infty\)  for all \(\BFx\in\calX\).
(Note that this does not mean the limit is necessarily zero.)
Not only is this result of interest in its own right, but also is crucial for proving the consistency of IKG policy under assumptions weaker than those imposed for previous related problems \citep{ScottFrazierPowell11}.
For example, the uniform convergence preserves the continuity of \(k_i^n(\BFx,\cdot)\) in the limit, a property that is crucial for the proof of \Cref{prop:step1}.
A more general version of \Cref{prop:unif_conv} is given in \citet[Proposition 2.9]{bect2019}, but we present a different proof built on RKHS theory in the Appendix.
\Cref{prop:unif_conv} sets a foundation for analysing the asymptotic behavior of Bayesian sequential sampling policies based on Gaussian processes with minimal assumptions.

Before we formally state \Cref{prop:step1,prop:step2}, the following definitions are required.
For each \(i\), let \(\eta_i^n\) denote the (random) number of times that a sample is taken from alternative \(i\) regardless of the value of \(\BFx\) up to the \(n\)-th sample, i.e.,
\[\eta_i^n\coloneqq \sum_{\ell=0}^{n-1}\ind_{\{a^\ell=i\}}.\]
Further, let \(\eta_i^\infty\coloneqq \lim_{n\to\infty} \eta_i^n\), which is well defined since it is a limit of a non-decreasing sequence of random variables.

\begin{proposition}\label{prop:step1}
    Fix \(i=1,\ldots,M\). If  \Cref{assump:compact,assump:stationary,assump:var} hold and \(\eta_i^\infty=\infty\) a.s., then for any \(\BFx\in\calX\), \(k_i^\infty(\BFx,\BFx)=0\) a.s. under the IKG policy.
\end{proposition}

\begin{proposition}\label{prop:step2}
    If  \Cref{assump:compact,assump:stationary,assump:var} hold, then \(\eta_i^\infty=\infty\) a.s.  for each \(i=1,\ldots,M\) under the IKG policy.
\end{proposition}

Part (i) of \Cref{theo:consistency} is an immediate consequence of  \Cref{prop:step1,prop:step2}.
The proofs of parts (ii) and (iii) of \Cref{theo:consistency} and \Cref{prop:step1,prop:step2} are all collected in the Appendix.

In practice, the IKG policy \eqref{eq:max-int-KG} can only be solved numerically, as discussed in the next section, in which case the obtained solution $(\tilde{a}^n,\tilde{\BFv}^n)$ is not exactly equal to the true solution $(a^n,\BFv^n)$.
Inspired by \cite{bect2019}, we consider the \emph{quasi}-IKG sampling policy, which chooses the sampling decision $(\tilde{a}^n,\tilde{\BFv}^n)$ such that
\begin{equation}\label{eq:max-int-KG-quasi}
\IKG^n(\tilde{a}^n,\tilde{\BFv}^n) \geq \IKG^n(a^n,\BFv^n) - \varepsilon_n,
\end{equation}
where $\{\varepsilon_n\}$ is a sequence of non-negative real numbers such that $\varepsilon_n\to 0$ as $n\to\infty$.
It is not difficult to see that such quasi-IKG policy is also consistent, as formally
stated in the following \Cref{theo:consistency-cost-quasi}, whose proof is collected in the Appendix.

\begin{theorem}\label{theo:consistency-cost-quasi}
    If \Cref{assump:stationary,assump:compact,assump:var} hold, then the quasi-IKG policy as defined in \cref{eq:max-int-KG-quasi} is  consistent.
\end{theorem}

\section{Stochastic Gradient Ascent}\label{sec:computation}

We now discuss computation of \cref{eq:max-int-KG} under \Cref{assump:stationary,assump:compact,assump:var}.
It primarily consists of two steps.
\begin{enumerate}[label=(\roman*)]
    \item \label{s1}
        For each \(i=1,\ldots,M\),  solve \(\max_{\BFx\in\calX} \IKG^n(i,\BFx)\) to find its maximizer, say \(\BFv^{n}_i\).
    \item \label{s2}
        Set \(a^n=\argmax_{1\leq i\leq M} \IKG^n(i, \BFv^{n}_i)\) and set  \(\BFv^n=\BFv^{n}_{a^n}\).
\end{enumerate}

Let \(\BFxi\) denote a \(\calX\)-valued random variable with density \(\gamma(\cdot)\), and
\begin{equation}\label{eq:function_h_i}
    h_i^n(\BFv, \BFx) \coloneqq
    \abs{\tilde\sigma_i^n(\BFv, \BFx)}\phi\qty(\abs{\frac{\Delta_i^n(\BFv)}{\tilde\sigma_i^n(\BFv, \BFx)}})
    -\abs{\Delta_i^n(\BFv)}\Phi\qty(-\abs{\frac{\Delta_i^n(\BFv)}{\tilde\sigma_i^n(\BFv, \BFx)}}).
\end{equation}
Then, we may rewrite \cref{eq:int-KG2} as
\begin{equation}\label{eq:IKG_expection_form}
    \IKG^n(i,\BFx) = [c_i(\BFx)]^{-1} \E[h_i^n(\BFxi, \BFx)],
\end{equation}
which suggests the following sample average approximation,
\begin{equation}\label{eq:sample_average}
    \widehat{\IKG}^n(i,\BFx) = \frac{1}{c_i(\BFx) J} \sum_{j=1}^{J} h_i^n(\BFxi_j, \BFx),
\end{equation}
where  \(\BFxi_j\)'s are independent copies of \(\BFxi\) and \(J\) is the sample size.
In particular, we will use \cref{eq:sample_average} in step (ii) above for computing \(a^n\) for given \(\BFv^n_i\)'s.
However, the sample average approximation method  can easily become computationally prohibitive when applied to  solve \(\max_{\BFx}\IKG^n(i,\BFx)\) in  step (i) if the domain \(\calX\) is high-dimensional. Hence, we consider instead the stochastic gradient ascent method to complete step (i).

\Cref{eq:IKG_expection_form} means that in step (i) above, we solve the stochastic optimization problem
\[\BFv^n_i \in \argmax_{\BFx\in\calX} \ [c_i(\BFx)]^{-1} \E \qty[h_i^n(\BFxi, \BFx)],\]
for each \(i=1,\ldots,M\).
If \(g_i^n(\BFxi, \BFx)\) is an unbiased estimator of \(\frac{\partial}{\partial \BFx} \{[c_i(\BFx)]^{-1} \E \qty[h_i^n(\BFxi, \BFx)]\}\), then
\(\BFv^n_i\) can be computed approximately using the stochastic gradient ascent (SGA) method; see~\cite{KushnerGeorge03} for a comprehensive treatment and~\cite{NewtonYousefianPasupathy18} for a recent survey on the subject.
Given an initial solution \(\BFx_1\in\calX\) and a maximum iteration limit \(K\), SGA iteratively computes
\begin{equation}\label{eq:SGA_ite}
    \BFx_{k+1} = \Pi_{\calX} \qty[\BFx_{k} + b_k g_i^n(\BFxi_k, \BFx_k) ],\quad k=1,\ldots,K,
\end{equation}
where \(\Pi_{\calX}:\Real^d\mapsto\calX\) denotes a projection mapping points outside \(\calX\) back to \(\calX\),\footnote{For example, one may set \(\Pi_{\calX}(\BFx)\) to be the point in \(\calX\) closest to \(\BFx\).}
and \(b_k\) is referred as the step size that satisfies \(\sum_{k=1}^\infty b_k=\infty\) and \(\sum_{k=1}^\infty b_k^2<\infty\).
In general, the choice of \(b_k\) is crucial for the practical performance of SGA, and it is commonly set as \(b_k = \alpha / k^{\beta}\) for some constants \(\alpha\) and \(\beta\).

Note that \(g_i^n(\BFxi, \BFx) = \frac{\partial}{\partial \BFx} [h_i^n(\BFxi, \BFx)/c_i(\BFx)]\) under mild regularity conditions~\citep{LEcuyer95}.
The explicit forms of $g_i^n(\BFxi, \BFx)$ for several common covariance functions are collected in the Appendix.
Besides, in the implementation of SGA algorithm, practical modifications such as mini batch and \emph{Polyak-Ruppert averaging}~\citep{PolyakJuditsky92} can be adopted to achieve better performance.
Detailed discussion is collected in the Appendix, together with other implementation issues of IKG policy.

\section{Numerical Experiments} \label{sec:numerical}

In this section, we evaluate the performance of the IKG policy via numerical experiments due to two reasons. First, the theoretical analysis, albeit establishing the consistency of the IKG policy in an large-sample asymptotic regime, does not provide a guarantee on the finite-sample performance of the policy. Second, the analysis has implicitly assumed that the sampling decisions of the IKG policy in \cref{eq:max-int-KG} can be computed exactly, while in practice it needs to be solved numerically via methods such as SGA that we have proposed.
Additional numerical experiments on other issues, including the computational cost comparison between SGA versus the sample average approximation and the effect of estimated $\lambda_i(\BFx)$, are collected in the Appendix.
All the numerical experiments are implemented in MATLAB and the source code is available at \url{https://github.com/shenhaihui/ikg}.

\subsection{Finite-Sample Performance}

The numerical experiments are conducted on synthetic problems, with the number of alternatives \(M=5\) and the dimensionality \(d=1,3,5,7\). For each \(i=1,\ldots,M\), the true performance of alternative \(i\) is the revised Griewank function,
\[
    \theta_i(\BFx) = \sum_{j=1}^d \frac{x_j^2}{4000} - 1.5^{d-1} \prod_{j=1}^d \cos\qty(\frac{x_j}{\sqrt{ij}}),\quad \BFx\in\calX=[0,10]^d.
\]
Further, we set sampling variance \(\lambda_i(\BFx) \equiv 0.01\), and take prior \(\mu_i^0(\BFx) = \mu^0(\BFx) \equiv 0\), and \(k_i^0(\BFx,\BFx') = k^0(\BFx,\BFx') = \exp\qty(- \frac{1}{d} \norm{\BFx-\BFx'}^2)\).
We set the cost function \(c_i(\BFx)\equiv 1\) for each \(i=1,\ldots,M\), but will investigate the impact of a different cost function later.

We consider two density functions for the covariates:
(1) uniform distribution on \(\calX\): \(\gamma(\BFx) = 1/\abs{\calX}\);
(2) multivariate normal distribution with mean \(\bm 0\) and covariance matrix \(4^2 \BFI\) truncated on \(\calX\): \(\gamma(\BFx) = \phi(\BFx;\bm 0, 4^2\BFI)/\int_{\calX} \phi(\BFv;\bm 0, 4^2\BFI) \dd{\BFv}\).
For convenience, we call the above specifications Problem 1 (P1) and Problem 2 (P2), respectively, depending on the choice of \(\gamma(\BFx)\).

The parameters involved in the SGA algorithm (see details in the Appendix) are given as follows: \(K=100d\), \(K_0 = K/4\), \(b_k = 200d / k^{0.7}\), \(m=20d\), and \(J=500d^2\). Moreover, the algorithm is started with a random initial solution. The performance of the IKG policy with respect to the sampling budget \(B\) is evaluated via the \emph{opportunity cost} (OC), that is, the integrated difference in performance between the best alternative and the  alternative chosen by the IKG policy upon exhausting the sampling budget.
\[\text{OC}(B)\coloneqq \E\qty[\int_{\calX}\qty(\theta_{i^*(\BFx)}(\BFx) - \theta_{\hat{i}^*(\BFx;\omega)}(\BFx))\gamma(\BFx)\dd{\BFx}],\]
where \(\hat{i}^*(\BFx;\omega) \in \argmax_{1\leq i\leq M}\mu^{N(B)}_i(\BFx;\omega)\) is the learned decision rule up to the budget \(B\) under the IKG policy, \(\omega\) denotes the samples taken under the policy, and the expectation is with respect to \(\omega\).
Clearly, \(\text{OC}(B)\to 0\) as  \(B\to\infty\), since the IKG policy is consistent. We estimate \(\text{OC}(B)\) via
\[
    \widehat{\text{OC}}(B) = \frac{1}{L} \sum_{l=1}^L \biggl[ \frac{1}{J'} \sum_{j=1}^{J'} \qty( \theta_{i^*(\BFx_j)}(\BFx_j) - \theta_{\hat{i}^*(\BFx_j;\omega_l)}(\BFx_j) ) \biggr],
\]
where \(L=30\) is the number of replications, \(\omega_l\) denotes the samples for replication \(l=1,\ldots,L\), and \(\{\BFx_1,\ldots,\BFx_{J'}\}\) is a random sample of the covariates generated from a given density function \(\gamma(\BFx)\) with \(J'=1000d^2\) for the purpose of evaluation.

\begin{figure}[!p]
    \begin{center}
        \caption{Estimated opportunity cost (vertical axis) as a function of the sampling budget (horizontal axis) for P1.} \label{fig:OC-P1}
            \includegraphics[width=0.36\textwidth]{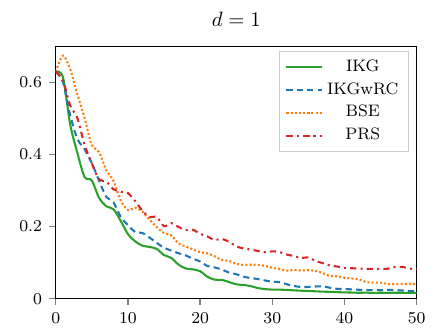} \hspace{0.05\textwidth}
            \includegraphics[width=0.36\textwidth]{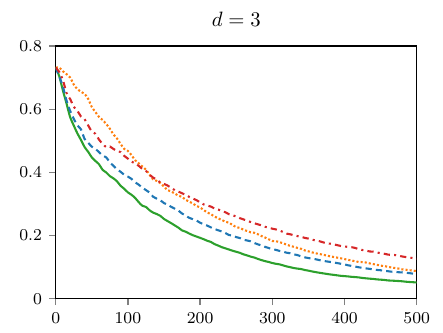}
            \includegraphics[width=0.36\textwidth]{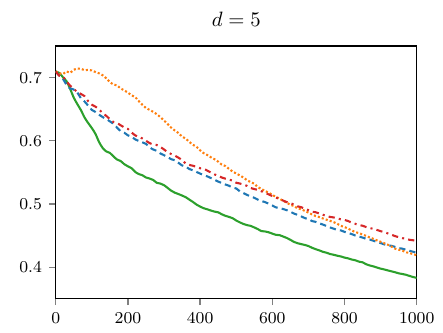} \hspace{0.05\textwidth}
            \includegraphics[width=0.36\textwidth]{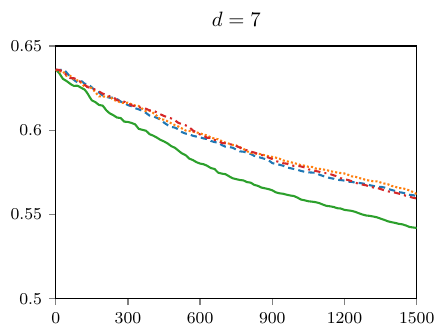}
    \end{center}
\end{figure}

\begin{figure}[!p]
    \begin{center}
        \caption{Estimated opportunity cost (vertical axis) as a function of the sampling budget (horizontal axis) for P2.} \label{fig:OC-P2}
            \includegraphics[width=0.36\textwidth]{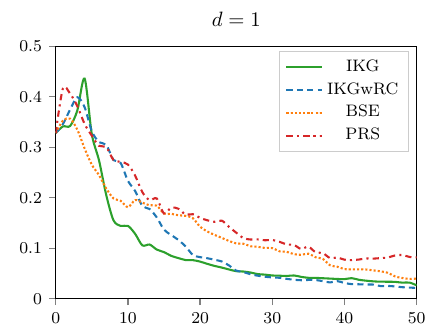} \hspace{0.05\textwidth}
            \includegraphics[width=0.36\textwidth]{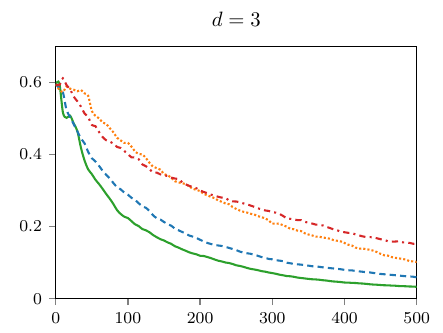}
            \includegraphics[width=0.36\textwidth]{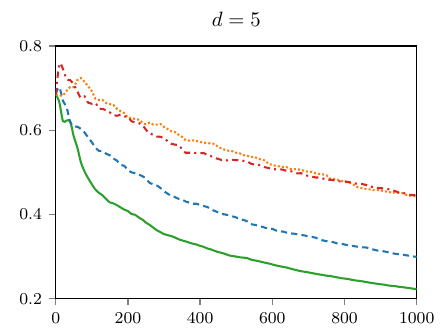} \hspace{0.05\textwidth}
            \includegraphics[width=0.36\textwidth]{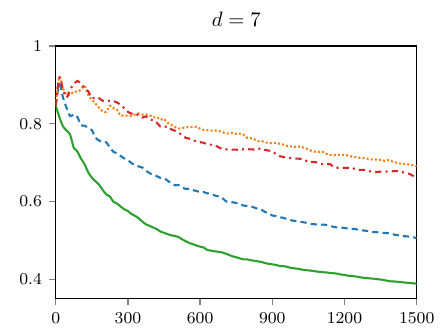}
    \end{center}
\end{figure}

We compare the IKG policy against three other polices:
\begin{itemize}
\item \emph{IKG with Random Covariates (IKGwRC).}
Recall that in the computation of IKG policy, random solution is used to initiate the SGA algorithm.
To check whether such random initialization is a main cause for the effectiveness of IKG, we consider the IKGwRC policy as follows.
Let $\BFx^n_1,\ldots,\BFx^n_M$ be the initial solutions for $M$ alternatives used in the SGA algorithm when computing \((\hat a^n, \hat \BFv^n)\), $n=0,1,\ldots$
Then the IKGwRC policy will sample at $(a^n,\BFv^n)$ given by
\[a^n= \argmax_{1\leq i\leq M} \log \widehat{\IKG}^n(i,\BFx^n_i)\qq{and}  \BFv^n= \BFx^n_{a^n},\]
where the same samples are used to compute $\widehat{\IKG}^n$ as in the IKG policy.

\item \emph{Binned Successive Elimination (BSE).}
The BSE policy is proposed by~\cite{perchet2013} for solving nonparametric MAB problems with covariates.
In their setting, values of the covaraites arrive randomly, and the policy only determines which alternative to select.
To implement BSE in our setting, we randomly generate $\BFv^n$ from uniform distribution on \(\calX\), and then apply the BSE policy to determine  $a^n$.
The BSE policy divides $\calX$ into $m^d$ parts, where $m$ is the number of uniformly divided regions on each coordinate.
For each problem, $m$ is tuned within the set $\{1,\ldots,10\}$, while other parameters follows the suggestion in~\cite{perchet2013}.

\item \emph{Pure Random Search (PRS).}
The PRS policy will sample at $(a^n,\BFv^n)$, where $a^n$ is randomly generated from the uniform distribution on $\{1,\ldots,M\}$ and $\BFv^n$ is generated from the uniform distribution on \(\calX\).
\end{itemize}
The performances of the four policies for problems P1 and P2 with $d=1,3,5,7$ are shown in \Cref{fig:OC-P1,fig:OC-P2}, respectively.  Several findings are made as follows.

First, the estimated opportunity cost in all the test problems exhibits a clear trend of convergence to zero.
This, from a practical point view, provides an assurance that the IKG policy in conjunction with the SGA algorithm indeed works as intended, that is, the uncertainty about the the performances of the competing alternatives will vanish eventually as the sampling budget grows. Second, the IKG policy can quickly reduce the opportunity cost when the sampling budget is relatively small, but the reduction appears to slow down as the sampling budget increases. This finding is consist with prior research on other KG-type policies such as \cite{FrazierPowellDayanik09}, \cite{FrazierPowell11}, and \cite{XieFrazierChick16}. Third, the learning task of identifying the best alternative becomes substantially more difficult when the dimensionality of the covariates is large. This can be seen from the growing sampling budget and the slowing reduction in the opportunity cost as \(d\) increases.

Overall, IKG outperforms the other three policies.
Specific comparisons are as follows.
First, IKG has better performance than IKGwRC, especially when the dimensionality is high, which indicates that the SGA algorithm in IKG for solving $\BFv_i^n$ (see \Cref{sec:computation}) indeed works well and has a significant effect in IKG.
Second, BSE has inferior performance than IKG, which may be caused by the fact that BSE only optimizes $a^n$ given randomly observed $\BFv^n$, while IKG optimizes both $a^n$ and $\BFv^n$ at the same time.
Third, PRS overall has the worst performance, which is not surprising since it does not utilize any information gained from previous sampling.
Note that PRS is a consistent policy, but the consistency does not guarantee any finite-sample performance.
This reflects the value of IKG -- it is not only provably consistent, but also takes advantage of information gained from previous samples to yield good finite-sample performance.

\subsection{Effect of Sampling Cost}

\begin{figure}[!ht]
    \begin{center}
        \caption{Estimated opportunity cost (vertical axis) as a function of the sampling budget (horizontal axis) for P3.} \label{fig:OC2}
            \includegraphics[width=0.36\textwidth]{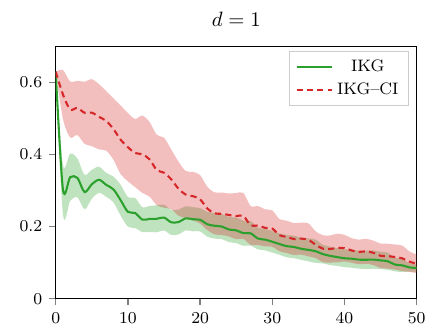} \hspace{0.05\textwidth}
            \includegraphics[width=0.36\textwidth]{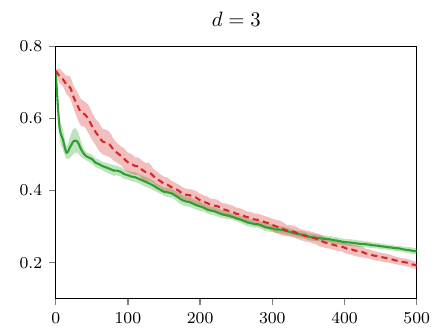}
            \includegraphics[width=0.36\textwidth]{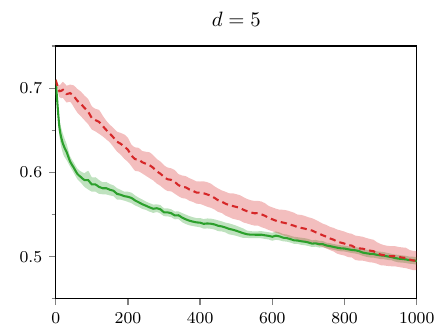} \hspace{0.05\textwidth}
            \includegraphics[width=0.36\textwidth]{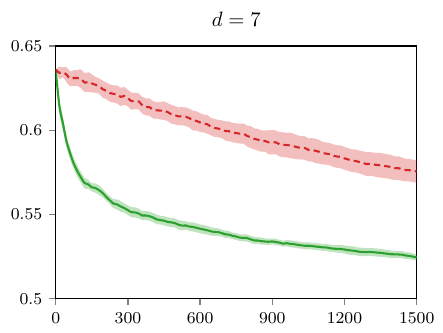}
    \end{center}
    \vspace{-10pt}
{\small \textsc{Note.} IKG--CI means sampling costs are ignored when implementing the IKG policy. The shaded regions represent the 99\% confidence intervals.}
\end{figure}

We are also interested in the effect of sampling costs on the IKG policy. In particular, we consider a different cost function other than the unit cost function: \(c_i(\BFx) = 2^{3-i} \bigl(1 + \norm{\BFx- \BFfive}^2/(10d) \bigr)\), where \(\BFfive\) is a \(d\times 1\) vector of all fives. We set \(\gamma(\BFx)\) to be the uniform density\footnote{Setting \(\gamma(\BFx)\) to be the truncated normal density leads to similar findings.} and call this specification Problem 3 (P3).
We compare two scenarios:
(i) the sampling cost is incorporated correctly;
and (ii) one ignores variations in the sampling cost at different locations and mistakenly uses the unit sampling cost when implementing the IKG policy (but the actual sampling consumption
follows \(c_i(\BFx)\)).
The comparison is illustrated in \Cref{fig:OC2}.

There are two observations. On one hand, despite the misspecification in the sampling cost function, the IKG policy is still consistent, with the associated opportunity cost converging to zero. This is not surprising, because using the unit sampling cost function, i.e., \(c_i(\BFx)\equiv 1\), is exactly the setup of \Cref{theo:consistency}. On the other hand, however, the finite-sample performance of the IKG policy indeed deteriorates as a result of the misspecification. Further, the deterioration appears to become more significant as the dimensionality of the covariates increases.

\section{Conclusions}\label{sec:conclusion}

In this paper, we study sequential sampling for the problem of selection with covariates which aims to identify the best alternative as a function of the covariates. Each sampling decision involves choosing an alternative and a value of the covariates, from the pair of which a sample will be taken. We design a sequential sampling policy via a nonparametric Bayesian approach. In particular, following the well-known KG design principle for simulation optimization, we develop the IKG policy that attempts to maximize the ``one-step'' integrated increment in the expected value of information per unit of sampling cost.

We prove the consistency of the IKG policy under minimal assumptions. Compared to prior work on asymptotic analysis of KG-type sampling policies, our assumptions are simpler and significantly more general, thanks to technical machinery that we develop based on RKHS theory. Nevertheless, to compute the sampling decisions of the IKG policy requires solving a multi-dimensional stochastic optimization problem. To that end, we develop a numerical algorithm based on the SGA method. Numerical experiments illustrate the finite-sample performance of the IKG policy and provide a practical assurance that the developed methodology works as intended.

\pdfbookmark[1]{References}{link-bib}
\bibliographystyle{chicago}
\bibliography{RSCov_KG}

\clearpage
\input{RSC-CIKG_appendix}

\end{document}

%% file: macros.tex
\DeclareMathOperator{\pr}{\mathbb{P}}
\DeclareMathOperator{\E}{\mathbb{E}}
\DeclareMathOperator{\Var}{Var}
\DeclareMathOperator{\Cov}{Cov}
\DeclareMathOperator{\Corr}{Corr}
\DeclareMathOperator{\ind}{\mathbb{I}}

\DeclareMathOperator*{\argmax}{argmax}
\def\Real{\mathbb{R}}

\def\BFA{\boldsymbol{A}}

\def\BFB{\boldsymbol{B}}

\def\BFI{\boldsymbol{I}}

\def\BFQ{\boldsymbol{Q}}

\def\BFV{\boldsymbol{V}}
\def\BFv{\boldsymbol{v}}

\def\BFX{\boldsymbol{X}}
\def\BFx{\boldsymbol{x}}

\def\BFy{\boldsymbol{y}}

\def\BFzero{\boldsymbol{0}}
\def\BFone{\boldsymbol{1}}

\def\BFfive{\boldsymbol{5}}

\def\BFalpha{\boldsymbol{\alpha}}

\def\BFdelta{\boldsymbol{\delta}}

\def\BFxi{\boldsymbol{\xi}}

\def\BFSigma{\boldsymbol{\Sigma}}

\def\calB{\mathcal{B}}
\def\calC{\mathcal{C}}

\def\calH{\mathcal{H}}

\def\calJ{\mathcal{J}}

\def\calN{\mathcal{N}}

\def\calS{\mathcal{S}}

\def\calX{\mathcal{X}}

\def\SFH{\mathsf{H}}

%% file: RSC-CIKG_appendix.tex
\appendix

\makeatletter
\renewcommand\@seccntformat[1]{}
\makeatother

\section{Appendix}

\subsection{A. Proof of Lemma \ref{lemma:calculating_IKG}}

Before proving \Cref{lemma:calculating_IKG}, we first establish the following \Cref{lemma:func_g}.

\begin{lemma}\label{lemma:func_g}
    Let \(g(s, t)\coloneqq t\phi(s/t)-s\Phi(-s/t)\), where \(\Phi\) is the standard normal distribution function and \(\phi\) is its density function. Then,
    \begin{enumerate}[label=(\roman*)]
        \item
              \(g(s,t)> 0\) for all \(s\geq 0\) and \(t> 0\);
        \item
              \(g(s,t)\) is strictly decreasing in \(s\in[0,\infty)\) and strictly increasing in \(t\in(0,\infty)\);
        \item
              \(g(s,t)\to 0\) as \(s\to\infty\) or as \(t\to 0\).
    \end{enumerate}
\end{lemma}

\begin{proof}[Proof of \Cref{lemma:func_g}.]
    Let \(h(u)\coloneqq \phi(u)-u\Phi(-u)\) for \(u\geq 0\), then \(g(s,t)=t h(s/t)\). Note that
    \[h'(u)=\phi'(u)+u\phi(-u) -\Phi(-u) = -u\phi(u) +u\phi(u) -\Phi(-u) = -\Phi(-u) <0, \]
    Hence, \(h(u)\) is strictly decreasing in \(u\in[0,\infty)\).
    Note that \(\lim_{u\to\infty} \phi(u) = 0\) and
    \[
        \lim_{u\to\infty} u\Phi(-u) = \lim_{u\to\infty} \frac{\Phi(-u)}{u^{-1}} = \lim_{u\to\infty} \frac{\phi(u)}{u^{-2}} = \lim_{u\to\infty} \frac{u^2}{\sqrt{2\pi} e^{u^2/2}} = 0,
    \]
    hence \(\lim_{u\to\infty} h(u) = 0\).
    Then we must have \(h(u)>0\) for all \(u\in[0,\infty)\), from which part (i) follows immediately.

    For part (ii), the strict decreasing monotonicity of \(g(s,t)\) in \(s\in[0,\infty)\) and the strict increasing monotonicity of \(g(s,t)\) in \(t\in(0,\infty)\) follow immediately from the strict decreasing monotonicity of \(h(u)\) in \(u\) and \(g(s,t)=t h(s/t)\).

    Part (iii) is due to that \(\lim_{u\to\infty} h(u) = 0\) and \(g(s,t)=t h(s/t)\).
\end{proof}

Now we are ready to prove \Cref{lemma:calculating_IKG}
\begin{proof}[Proof of \Cref{lemma:calculating_IKG}.]

    By \cref{eq:mean_update,eq:sigma},
    \begin{align}
        f(i,\BFx,\BFv) & \coloneqq \E\qty[\max_{1\leq a\leq M} \mu^{n+1}_a(\BFv)\,\Big|\,\mathscr{F}^n, a^n=i, \BFv^n=\BFx] \nonumber                              \\[0.5ex]
                       & = \E\qty[\max_{1\leq a\leq M}\qty(\mu_a^n(\BFv)+\sigma_a^n(\BFv, \BFx)Z^{n+1})] \nonumber                                          \\[0.5ex]
                       & = \E\qty[\max\qty{\mu_i^n(\BFv) +\tilde\sigma_i^n(\BFv, \BFx)Z^{n+1}, \;\max_{a\neq i}\mu_a^n(\BFv)}] \nonumber                  \\[0.5ex]
                       & = \E\qty[\max\qty{\mu_i^n(\BFv) + |\tilde\sigma_i^n(\BFv, \BFx)| Z^{n+1}, \;\max_{a\neq i}\mu_a^n(\BFv)}]. \label{eq:f_function}
    \end{align}
    For notational simplicity, let \(\alpha \coloneqq \mu_i^n(\BFv)\), \(\beta \coloneqq \tilde\sigma_i^n(\BFv, \BFx)\), \(\gamma \coloneqq \max_{a\neq i}\mu_a^n(\BFv)\), and \(\delta \coloneqq \alpha-\gamma\).

    If \(\beta\neq0\), then
    \begin{align*}
        f(i,\BFx,\BFv) = &  \int_{-\infty}^{-\delta/|\beta|} \gamma \phi(z)\dd{z}
        + \int_{-\delta/|\beta|}^\infty (\alpha+|\beta| z)\phi(z)\dd{z} \\
        = & \gamma \Phi(-\delta /|\beta| ) + \alpha[1-\Phi(-\delta/|\beta|)] + \beta \phi(-\delta/|\beta|),
    \end{align*}
    where the second equality follows from the identity \(\int_t^\infty z\phi(z)\dd{z}=\phi(t)\) for all \(t\in\Real\).
    Next, we calculate the integrand in \cref{eq:int-KG}.
    By noting that \(\phi(z)=\phi(-z)\) and \(\Phi(-z)=1-\Phi(z)\) for all \(z\in\Real\),
    \begin{align}
        f(i,\BFx,\BFv)-\max_{1\leq a\leq M}\mu_a^n(\BFv)  = f(i,\BFx,\BFv)- \max(\alpha,\gamma) & =
        \begin{cases}
            |\beta|\phi(\delta/|\beta|)- \delta\Phi(-\delta/|\beta|), & \text{if } \alpha \geq \gamma \\
            |\beta|\phi(\delta/|\beta|)+\delta\Phi(\delta/|\beta|),   & \text{if } \alpha < \gamma
        \end{cases}
        \nonumber  \\
                                                                                                & = |\beta|\phi(|\delta/\beta|)-|\delta|\Phi(-|\delta/\beta|).   \label{eq:IKG_beta_nonzero}
    \end{align}
    If \(\beta=0\), then it is straightforward from \cref{eq:f_function} to see that \(f(i,\BFx,\BFv)-\max_{1\leq a\leq M}\mu_a^n(\BFv) =0\). On the other hand, by \Cref{lemma:func_g} (iii), we can set the right-hand side of \cref{eq:IKG_beta_nonzero} to be zero for \(\beta=0\). Hence, \cref{eq:IKG_beta_nonzero} holds for \(\beta=0\) as well.

    Replacing the integrand in \cref{eq:int-KG} with \cref{eq:IKG_beta_nonzero} yields the expression of \(\IKG^n(i,\BFx)\) in \Cref{lemma:calculating_IKG}.
\end{proof}

\subsection{B. Proof of Proposition \ref{prop:unif_conv}}\label{sec:RKHS}

To simplify notation, in this subsection we assume  \(M=1\) and suppress the subscript \(i\) unless otherwise specified, but the results can be generalized to the case of \(M>1\) without essential difficulty. In particular, we use \(\kappa\) to denote a generic covariance function, \(k^0\) the prior covariance function of a Gaussian process, and \(k^n\) the posterior covariance function.
We will collect below several basic results on reproducing kernel Hilbert space (RKHS) and refer to~\cite{BerlinetThomas-Agnan04} for an extensive treatment on the subject.

\begin{definition}\label{def:RKHS}
    Let \(\calX\) be a nonempty set and \(\kappa\) be a covariance function on \(\calX\).
    A Hilbert space \(\calH_\kappa\) of functions on \(\calX\) equipped with an inner-product \(\langle\cdot,\cdot\rangle_{\calH_\kappa}\) is called a RKHS with reproducing kernel \(\kappa\), if (i) \(\kappa(\BFx, \cdot)\in\calH_\kappa\) for all \(\BFx\in\calX\), and (ii) \(f(\BFx)=\langle f,\kappa(\BFx, \cdot)\rangle_{\calH_\kappa}\) for all \(\BFx\in\calX\) and \(f\in\calH_\kappa\).
    Furthermore, the norm of \(\calH_\kappa\) is induced by the inner-product, i.e., \(\norm{f}_{\calH_\kappa}^2 =\langle f,f\rangle_{\calH_\kappa}\) for all \(f\in\calH_\kappa\).
\end{definition}
\begin{remark}
    In \Cref{def:RKHS}, for a fixed \(\BFx\), \(\kappa(\BFx, \cdot)\) is understood as a function mapping \(\calX\) to \(\Real\) such that \(\BFy\mapsto k(\BFx,\BFy)\) for \(\BFy\in\calX\). Moreover, condition (ii) is  called the \emph{reproducing property}.
    In particular, it implies that \(\kappa(\BFx,\BFx')=\langle \kappa(\BFx, \cdot), \kappa(\BFx', \cdot)\rangle_{\calH_\kappa}\) and \(\kappa(\BFx,\BFx)=\norm{\kappa(\BFx, \cdot)}_{\calH_\kappa}^2\) for all \(\BFx,\BFx'\in\calX\).
\end{remark}
\begin{remark}
    By Moore-Aronszajn theorem~\cite[Theorem 3]{BerlinetThomas-Agnan04}, for each covariance function \(\kappa\) there exists a unique RKHS \(\calH_\kappa\) for which \(\kappa\) is its reproducing kernel. Specifically,
    \[\calH_\kappa=\qty{f=\sum_{i=1}^\infty c_i \kappa(\BFx_i, \cdot): c_i\in\Real, \BFx_i\in\calX,\,i=1,2,\ldots, \text{ such that } \norm{f}^2_{\calH_\kappa}<\infty},\]
    where \(\norm{f}^2_{\calH_\kappa}\coloneqq \sum_{i,j=1}^\infty c_ic_j \kappa(\BFx_i,\BFx_j)\).
    Moreover, the inner-product is defined by
    \[\langle f,g\rangle_{\calH_\kappa} = \sum_{i,j=1}^\infty a_ib_j \kappa(\BFx_i,\BFx'_j),\]
    for any \(f=\sum_{i=1}^\infty a_i \kappa(\BFx_i, \cdot)\in\calH_\kappa\) and \(g=\sum_{j=1}^\infty b_j \kappa(\BFx_j', \cdot)\in\calH_\kappa\).
\end{remark}

The following lemma asserts that convergence in norm in a RKHS implies \emph{uniform} pointwise convergence, provided that the covariance function \(\kappa\) is stationary.

\begin{lemma}\label{lemma:unif_conv}
    Let \(\calX\) be a nonempty set and \(\kappa\) be a covariance function on \(\calX\).
    Suppose that  a sequence of functions \(\{f_n\in\calH_\kappa:n=1,2,\ldots\}\) converges in norm \(\norm{\cdot}_{\calH_\kappa}\) as $n \to \infty$.
    Then the limit, denoted by $f$, is in \(\calH_\kappa\).
    Moreover, if \(\kappa\) is stationary, then \(f_n(\BFx)\to f(\BFx)\) as \(n\to\infty\) uniformly in \(\BFx\in\calX\).
\end{lemma}

\begin{proof}[Proof of \Cref{lemma:unif_conv}.]
First of all, \(f \in \calH_\kappa\) is guaranteed as a Hilbert space is a complete metric space.
A basic property of RKHS is that convergence in norm implies pointwise convergence to the same limit; see, e.g., Corollary 1 of~\citet[page 10]{BerlinetThomas-Agnan04}. Namely, \(f_n(\BFx)\to f(\BFx)\) as \(n\to\infty\) for all \(\BFx\in\calX\).

    To show the pointwise convergence is uniform, note that since \(\kappa\) is stationary, there exists a function \(\varphi:\Real^d \mapsto\Real\) such that \(\kappa(\BFx,\BFx')=\varphi(\BFx-\BFx')\).
    Hence, \(\norm{\kappa(\BFx, \cdot)}_{\calH_\kappa}^2=\kappa(\BFx,\BFx)=\varphi(\BFzero)\).
    It follows that
    \begin{align}
        \abs{f_{n+m}(\BFx)-f_n(\BFx)} & = \abs{\langle f_{n+m}-f_n, \kappa(\BFx,\cdot)\rangle_{\calH_\kappa}} \nonumber\\
        & \leq \norm{f_{n+m}-f_n}_{\calH_\kappa}\norm{\kappa(\BFx, \cdot)}_{\calH_\kappa}=\norm{f_{n+m}-f_n}_{\calH_\kappa} \sqrt{\varphi(\BFzero)}, \label{eq:Cauchy_seq}
    \end{align}
    for all \(n\) and \(m\), where the first equality follows from the reproducing property.

    Since a Hilbert space is a complete metric space, the \(\norm{\cdot}_{\calH_\kappa}\)-converging sequence \(\{f_n\}\) is a Cauchy sequence in \(\calH_\kappa\), meaning that \(\norm{f_{n+m}-f_n}_{\calH_\kappa}\to 0\) as \(n\to\infty\) for all \(m\). Since this convergence to zero is independent of \(\BFx\), it follows  from \cref{eq:Cauchy_seq} that \(\{f_n\}\) is a uniform Cauchy sequence of functions, thereby  converging to \(f\) uniformly in \(\BFx\in\calX\).
\end{proof}

In the light of \Cref{lemma:unif_conv}, in order to establish the uniform convergence of \(k^n(\BFx,\BFx')\) as a function of \(\BFx'\), it suffices to prove the norm convergence of \(k^n(\BFx,\cdot)\) in the RKHS induced by \(k^0\).
We first establish this result for a more general case in the following \Cref{lemma:Cauchy_seq}, where \(k^0\) is not required to be stationary.

\begin{lemma}\label{lemma:Cauchy_seq}
    Let \(\calH_{k^0}\) be the RKHS induced by \(k^0\). If \(k^0(\BFx,\BFx)>0\) for all \(\BFx\in\calX\), then for any \(\BFx\in\calX\), \(k^n(\BFx,\cdot)\) converges in norm \(\norm{\cdot}_{\calH_{k^0}}\) as \(n\to\infty\).
\end{lemma}

\begin{proof}[Proof of \Cref{lemma:Cauchy_seq}.]

    Fix \(\BFx\in\calX\).
    The fact that \(k^n(\BFx,\cdot) \in \calH_{k^0}\) is due to \cref{eq:cov_n}.
    It follows from \cref{eq:decreasing_var} that \(\{k^n(\BFx,\BFx):n\geq 1\}\) form a non-increasing sequence bounded below by zero.
    The monotone convergence theorem implies that \(k^n(\BFx,\BFx)\) converges as \(n\to\infty\). Hence, for all \(m\geq 1\),
    \begin{equation}\label{eq:Cauchy_pointwise}
        \lim_{n\to\infty} \abs{k^{n+m}(\BFx,\BFx)-k^n(\BFx,\BFx)}=0.
    \end{equation}

    Let \(\BFV^n\coloneqq \{\BFv^\ell:\ell=0,\ldots,n-1\}\) and \(\BFV_n^{n+m}\coloneqq \{\BFv^\ell: \ell=n,\ldots,n+m-1\}\). Then, by \cref{eq:cov_n},
    \begin{equation}\label{eq:update_from_n}
        k^{n+m}(\BFx, \cdot) - k^n(\BFx, \cdot) = - k^n(\BFx,\BFV_n^{n+m})[k^n(\BFV_n^{n+m},\BFV_n^{n+m})+\lambda(\BFV_n^{n+m})]^{-1}k^n(\BFV_n^{n+m},\cdot).
    \end{equation}
    For notational simplicity, let \(\BFSigma_n^{n+m}\coloneqq k^n(\BFV_n^{n+m},\BFV_n^{n+m})+\lambda(\BFV_n^{n+m})\). Then,
    \begin{align}
          & \norm{k^{n+m}(\BFx,\cdot)-k^n(\BFx,\cdot)}_{\calH_{k^0}}^2  \nonumber                                                                                                                                         \\[0.5ex]
        = {} & \,
        \langle k^n(\BFx,\BFV_n^{n+m}) [\BFSigma_n^{n+m}]^{-1} k^n(\BFV_n^{n+m},\cdot), k^n(\BFx,\BFV_n^{n+m}) [\BFSigma_n^{n+m}]^{-1} k^n(\BFV_n^{n+m},\cdot) \rangle_{\calH_{k^0}} \nonumber                         \\[0.5ex]
        = {} & \, k^n(\BFx,\BFV_n^{n+m}) [\BFSigma_n^{n+m}]^{-1} \langle k^n(\BFV_n^{n+m},\cdot),k^n(\BFV_n^{n+m},\cdot) \rangle_{\calH_{k^0}} [\BFSigma_n^{n+m}]^{-1} k^n(\BFV_n^{n+m},\BFx), \label{eq:norm_calcuation}
    \end{align}
    where
    \[\langle k^n(\BFV_n^{n+m},\cdot),k^n(\BFV_n^{n+m},\cdot) \rangle_{\calH_{k^0}} = \qty(\langle k^n(\BFv^{n+i},\cdot), k^n(\BFv^{n+j},\cdot)  \rangle_{\calH_{k^0}})_{0\leq i,j\leq m-1}.\]
    Moreover, note that by \cref{eq:cov_n},
    \begin{equation}\label{eq:complicated}
        k^n(\BFV_n^{n+m},\cdot) = k^0(\BFV_n^{n+m},\cdot)-k^0(\BFV_n^{n+m}, \BFV^n)[k^0(\BFV^n,\BFV^n)+\lambda(\BFV^n)]^{-1}k^0(\BFV^n,\cdot).
    \end{equation}
    Let \(\BFSigma^n\coloneqq k^0(\BFV^n,\BFV^n)+\lambda(\BFV^n)\). Then, it follows from \cref{eq:complicated} and the reproducing property that
    \begin{align}
          & \langle k^n(\BFV_n^{n+m},\cdot),k^n(\BFV_n^{n+m},\cdot) \rangle_{\calH_{k^0}} \nonumber                                                                                        \\[0.5ex]
        = {} & \, k^0(\BFV_n^{n+m},\BFV_n^{n+m})-2 k^0(\BFV_n^{n+m}, \BFV^n)[\BFSigma_n]^{-1}k^0(\BFV^n,\BFV_n^{n+m}) \nonumber                                                               \\[0.5ex]
          & +\, k^0(\BFV_n^{n+m}, \BFV^n)[\BFSigma_n]^{-1} k^0(\BFV^n,\BFV^n)[\BFSigma_n]^{-1} k^0(\BFV^n,\BFV_n^{n+m})\nonumber                                                           \\[0.5ex]
        = {} & \, k^0(\BFV_n^{n+m},\BFV_n^{n+m})- k^0(\BFV_n^{n+m}, \BFV^n)[\BFSigma_n]^{-1}k^0(\BFV^n,\BFV_n^{n+m}) \nonumber                                                                \\[0.5ex]
          & -\, k^0(\BFV_n^{n+m}, \BFV^n)[\BFSigma_n]^{-1} \{\BFI- k^0(\BFV^n,\BFV^n)[\BFSigma_n]^{-1}\} k^0(\BFV^n,\BFV_n^{n+m}) \nonumber                                                \\[0.5ex]
        = {} & \, k^n(\BFV_n^{n+m},\BFV_n^{n+m}) - k^0(\BFV_n^{n+m}, \BFV^n)[\BFSigma_n]^{-1} \{\BFI- k^0(\BFV^n,\BFV^n)[\BFSigma_n]^{-1}\} k^0(\BFV^n,\BFV_n^{n+m}), \label{eq:complicated2}
    \end{align}
    where \(\BFI\) denotes the identity matrix of a compatible size.
    Furthermore, note that
    \begin{align}
        \BFI- k^0(\BFV^n,\BFV^n)[\BFSigma_n]^{-1} & = \BFI-[k^0(\BFV^n,\BFV^n)+\lambda(\BFV^n)-\lambda(\BFV^n)][k^0(\BFV^n,\BFV^n)+\lambda(\BFV^n)]^{-1} \nonumber \\[0.5ex]
                                                  & = \BFI- \BFI+\lambda(\BFV^n) [k^0(\BFV^n,\BFV^n)+\lambda(\BFV^n)]^{-1} \nonumber                               \\[0.5ex]
                                                  & = \lambda(\BFV^n)[\BFSigma^n]^{-1}.\label{eq:complicated3}
    \end{align}
    We now combine \cref{eq:complicated2,eq:complicated3} to have
    \begin{align*}
          & \langle k^n(\BFV_n^{n+m},\cdot),k^n(\BFV_n^{n+m},\cdot) \rangle_{\calH_{k^0}} \\[0.5ex]
        ={} & \,
        k^n(\BFV_n^{n+m},\BFV_n^{n+m}) -  k^0(\BFV_n^{n+m}, \BFV^n)[\BFSigma_n]^{-1} \lambda(\BFV^n)[\BFSigma^n]^{-1} k^0(\BFV^n,\BFV_n^{n+m}),
    \end{align*}
    which is the difference between two positive semi-definite matrices.
    Therefore, by \cref{eq:norm_calcuation},
    \begin{align*}
        \norm{k^{n+m}(\BFx,\cdot)-k^n(\BFx,\cdot)}_{\calH_{k^0}}^2
         & \leq
        k^n(\BFx,\BFV_n^{n+m}) [\BFSigma_n^{n+m}]^{-1} k^n(\BFV_n^{n+m},\BFV_n^{n+m}) [\BFSigma_n^{n+m}]^{-1} k^n(\BFV_n^{n+m},\BFx) \\[0.5ex]
         & \leq k^n(\BFx,\BFV_n^{n+m}) [\BFSigma_n^{n+m}]^{-1} \BFSigma_n^{n+m} [\BFSigma_n^{n+m}]^{-1} k^n(\BFV_n^{n+m},\BFx)       \\[0.5ex]
         & = k^n(\BFx,\BFx) - k^{n+m}(\BFx,\BFx),
    \end{align*}
    where the second inequality follows from the definition of \(\BFSigma_n^{n+m}\) and the equality follows from \cref{eq:update_from_n}. Then, we apply \cref{eq:Cauchy_pointwise} to conclude that
    \(\norm{k^{n+m}(\BFx,\cdot)-k^n(\BFx,\cdot)}_{\calH_{k^0}}\to 0\) as \(n\to\infty\)
    for all \(m\geq 1\). Therefore, \(k^n(\BFx,\cdot)\) converges in norm \(\norm{\cdot}_{\calH_{k^0}}\) as \(n\to\infty\).
\end{proof}

With \Cref{lemma:unif_conv,lemma:Cauchy_seq}, we are ready to prove \Cref{prop:unif_conv}.

\begin{proof}[Proof of \Cref{prop:unif_conv}.]
Fix \(i=1,\ldots,M\).
Since \(k_i^0\) is stationary, \(k_i^0(\BFx,\BFx)>0\) for all \(\BFx\in\calX\).
Then by \Cref{lemma:Cauchy_seq}, for any \(\BFx\in\calX\), \(k_i^n(\BFx,\cdot)\) converges in norm \(\norm{\cdot}_{\calH_{k^0}}\) as \(n\to\infty\).
Then by \Cref{lemma:unif_conv}, for the \(\norm{\cdot}_{\calH_{k^0}}\)-converging limit \(k_i^\infty(\BFx,\cdot)\),
\(k_i^n(\BFx,\BFx') \to k_i^\infty(\BFx,\BFx')\) uniformly in \(\BFx'\in\calX\) as \(n\to\infty\).
\end{proof}

\subsection{C. Proof of Proposition \ref{prop:step1}}

Notice that if \(\eta_i^\infty=\infty\) under a sampling policy \(\pi\), then due to the compactness of \(\calX\), \(\{\BFv^n:a^n=i,n=0,1,\ldots\}\) (i.e., the sampling locations associated with alternative \(i\)) must have an accumulation point \(\BFx_i^{\acc}\in\calX\).
Namely, there exists a subsequence of \(\{n:a^n=i, n=0,1,\ldots\}\), say \(\{\ell_{i,n}\}_{n=0}^\infty\), such that \(\ell_{i,n}\to\infty\) and \(\BFv^{\ell_{i,n}}\to \BFx_i^{\acc}\) as \(n\to\infty\).
For any \(\epsilon>0\), let \(\calB(\BFx_i^{\acc}, \epsilon)\coloneqq \{\BFx:\norm{\BFx_i^{\acc}-\BFx}\leq \epsilon\}\) be the closed ball centered at \(\BFx_i^{\acc}\) with radius \(\epsilon\).
Let \(\Var^{\pi, n}[\cdot]\) denote the posterior variance conditioned on \(\mathscr{F}^n\) that is induced by \(\pi\).

The proof of \Cref{prop:step1} is preceded by four technical results, i.e., \Cref{lemma:bound_acc_1,lemma:bound_acc_2,lemma:liminf,lemma:KG_bound}.
In \Cref{lemma:bound_acc_1,lemma:bound_acc_2}, we establish an upper bound on \(\Var^{\pi,n}[\theta_i(\BFx)]\) for \(\BFx\in \calB(\BFx_i^{\acc}, \epsilon)\).
This result does not rely on the IKG policy \emph{per se}, but is implied by the existence of the accumulation point \(\BFx_i^{\acc}\) instead.
In particular, the upper bound which depends on \(\epsilon\) can be made arbitrarily small as \(\epsilon\to 0\).
This basically means that in the light of an unlimited number of samples of alternative \(i\) that are taken in proximity to \(\BFx_i^{\acc}\), the uncertainty about \(\theta_i(\BFx_i^{\acc})\) will eventually be eliminated,
thanks to the correlation between \(\theta_i(\BFx_i^{\acc})\) and \(\theta_i(\BFx)\) for \(\BFx\) near \(\BFx_i^{\acc}\).

\Cref{lemma:KG_bound} asserts that \(\IKG^n(i,\BFx)\) is bounded by a multiple of the posterior standard deviation of \(\theta_i(\BFx)\).
This implies that when the posterior variance approaches to zero, the IKG factor does too.

Following  the last three lemmas, \Cref{lemma:liminf} asserts that the limit inferior of the IKG factor is zero.
The reasoning is as follows.
By \Cref{lemma:bound_acc_1,lemma:bound_acc_2}, the posterior variance at those sampling locations that fall inside \(\calB(\BFx_i^{\acc}, \epsilon)\) is small. Then, by \Cref{lemma:KG_bound}, the IKG factor at these locations are also small, so does the limit superior.
Since the sampling locations inside \(\calB(\BFx_i^{\acc}, \epsilon)\) is a subsequence of the entire sampling locations, the limit inferior of the IKG factor over all sampling locations is even smaller.

At last, \Cref{prop:step1} is proved by contradiction --- if there is a location that the posterior variance does not approach zero, then the limit inferior of the IKG factor at the same location must be positive as well.

\begin{lemma}\label{lemma:bound_acc_1}
    Fix \(i=1,\ldots,M\), \(n\geq 1\), and a compact set \(\calS\subseteq\calX\). Suppose that the sampling decisions satisfy \(a^0=\cdots=a^{n-1}=i\) and \(\BFv^0,\ldots,\BFv^{n-1}\in\calS\).     If \Cref{assump:stationary,assump:var,assump:compact} hold, then for all \(\BFx\in\calS\),
    \[\Var^n[\theta_i(\BFx)] \leq \tau_i^2 - \frac{n \min_{\BFx' \in \calS} \qty[k_i^0(\BFx,\BFx')]^2}{n \tau_i^2+\lambda_i^{\max}}, \]
    where \(\lambda_i^{\max}\coloneqq  \max_{\BFx\in \calX} \lambda_i(\BFx)\in(0,\infty)\).
\end{lemma}

\begin{proof}[Proof of \Cref{lemma:bound_acc_1}.]
    Fix \(\BFx\in\calS\).
    First note that \(\lambda_i^{\max}\) is well defined under \Cref{assump:var,assump:compact}.
    Let \(\BFV_i^n\) be the set of the locations of the samples taken from \(\theta_i\) up to time \(n\).
   Under \Cref{assump:stationary}, \cref{eq:cov_n} reads
    \[\Var^n[\theta_i(\BFx)] = \tau_i^2 - k_i^0(\BFx,\BFV_i^n)[k_i^0(\BFV_i^n,\BFV_i^n)+\lambda_i(\BFV_i^n)]^{-1}k_i^0(\BFV_i^n,\BFx),\]
    where \(\BFV_i^n=\{\BFv^0,\ldots,\BFv^{n-1}\}\) due to the assumption that \(a^0=\cdots=a^{n-1}=i\).

For notational simplicity, let \(\BFA \coloneqq k_i^0(\BFV_i^n,\BFV_i^n)+ \lambda_i(\BFV_i^n)\) and \(\BFB \coloneqq k_i^0(\BFV_i^n,\BFV_i^n) +  \lambda_i^{\max}\BFI\).
Note that \(\BFB-\BFA = \lambda_i^{\max}\BFI - \lambda_i(\BFV_i^n)\) is a diagonal matrix with nonnegative elements, so it is positive semi-definite.
Since \(\BFA\) and \(\BFB\) are both positive definite, by \citet[Corollary 7.7.4]{HornJohnson12}, \(\BFA^{-1}-\BFB^{-1}\) is positive semi-definite.
Therefore,
\begin{align}
&  k_i^0(\BFx,\BFV_i^n)[k_i^0(\BFV_i^n,\BFV_i^n)+\lambda_i(\BFV_i^n)]^{-1}k_i^0(\BFV_i^n,\BFx)  - k_i^0(\BFx,\BFV_i^n)[k_i^0(\BFV_i^n,\BFV_i^n)+\lambda_i^{\max}\BFI]^{-1}k_i^0(\BFV_i^n,\BFx) \nonumber\\
= {} & k_i^0(\BFx,\BFV_i^n) (\BFA^{-1}-\BFB^{-1}) k_i^0(\BFV_i^n,\BFx)
\geq 0. \label{eq:PSD}
\end{align}

    It then follows from \cref{eq:cov_n,eq:PSD} that
    \[\Var^n[\theta_i(\BFx)]  \leq \tau_i^2 - k_i^0(\BFx,\BFV_i^n)[k_i^0(\BFV_i^n,\BFV_i^n)+\lambda_i^{\max}\BFI]^{-1}k_i^0(\BFV_i^n,\BFx).\]
    Thus, it suffices to prove that
    \begin{equation}\label{eq:var_up_bound2}
    f(\BFv^0,\ldots,\BFv^{n-1}) \coloneqq k_i^0(\BFx,\BFV_i^n)[k_i^0(\BFV_i^n,\BFV_i^n)+ \lambda_i^{\max} \BFI]^{-1}k_i^0(\BFV_i^n,\BFx) \geq \frac{n \min_{\BFx'\in\calS} \qty[k_i^0(\BFx,\BFx')]^2}{n \tau_i^2+\lambda_i^{\max}},
    \end{equation}
    for all \(\BFv^0,\ldots,\BFv^{n-1}\in\calS\).

    Since \(k_i^0(\BFV_i^n,\BFV_i^n)\) is symmetric, we can always write \(k_i^0(\BFV_i^n,\BFV_i^n) = \BFQ \mathrm{diag}\{\alpha_1,\ldots,\alpha_n\} \BFQ^\intercal\), where \(\alpha_1\geq \alpha_2 \geq \cdots \geq \alpha_n\geq 0\) are the eigenvalues of \(k_i^0(\BFV_i^n,\BFV_i^n)\), and \(\BFQ\) is an orthogonal matrix, i.e., \(\BFQ \BFQ^\intercal = \BFI\).
Therefore,
    \[
    k_i^0(\BFV_i^n,\BFV_i^n)+ \lambda_i^{\max} \BFI = \BFQ \mathrm{diag} \big\{ (\alpha_1+\lambda_i^{\max}), \ldots, (\alpha_n+\lambda_i^{\max}) \big\} \BFQ^\intercal.
    \]
and
    \[
    [k_i^0(\BFV_i^n,\BFV_i^n)+ \lambda_i^{\max} \BFI]^{-1} = \BFQ \mathrm{diag} \big\{ (\alpha_1+\lambda_i^{\max})^{-1}, \ldots, (\alpha_n+\lambda_i^{\max})^{-1} \big\} \BFQ^\intercal.
    \]
    If we let \(\beta_j\) be the \(j\)-th element of the row vector \(k_i^0(\BFx,\BFV_i^n) \BFQ\), i.e., \(k_i^0(\BFx,\BFV_i^n) \BFQ = [\beta_1,\ldots,\beta_n]\), then
    \[f(\BFv^0,\ldots,\BFv^{n-1}) = \frac{\beta_1^2}{\alpha_1+\lambda_i^{\max}} + \cdots + \frac{\beta_n^2}{\alpha_n+\lambda_i^{\max}}.\]
    Here, \(\alpha_j\) and \(\beta_j\) clearly both depend on \(\BFv^0,\ldots,\BFv^{n-1}\), for \(j=1,\ldots,n\). Moreover, they satisfy the following two conditions. First,  \(\sum_{j=1}^n \alpha_j = \trace(k_i^0(\BFV_i^n,\BFV_i^n)) = n \tau_i^2\), where the first equality is a straightforward fact that the trace of a matrix equals the sum of its eigenvalues,
    and the second equality is from \Cref{assump:stationary}. Second,
    \[\sum_{j=1}^n \beta_j^2 = k_i^0(\BFx,\BFV_i^n) \BFQ \BFQ^\intercal k_i^0(\BFV_i^n,\BFx)
     = k_i^0(\BFx,\BFV_i^n) k_i^0(\BFV_i^n,\BFx) = \sum_{\ell=0}^{n-1} [k_i^0(\BFx,\BFv^\ell)]^2 \geq n \min_{\BFx'\in\calS}[k_i^0(\BFx,\BFx')]^2.\]

     If we define \(g:\Real^{2n}\mapsto \Real\) as follows
     \[g(a_1,\ldots,a_n,b_1,\ldots,b_n) \coloneqq \frac{b_1}{a_1+\lambda_i^{\max}} + \cdots + \frac{b_n}{a_n+\lambda_i^{\max}},\]
     then \(f(\BFv^0,\ldots,\BFv^{n-1}) = g(\alpha_1,\ldots,\alpha_n,\beta_1^2,\ldots,\beta_n^2)\).
It follows that
    \begin{equation}\label{eq:min-p-extend}
        \min_{\BFv^0,\ldots,\BFv^{n-1} \in \calS} f(\BFv^0,\ldots,\BFv^{n-1}) \geq  \min_{\substack{(a_1,\ldots,a_n)\in\calC_1 \\ (b_1,\ldots,b_n)\in \calC_2}} g(a_1,\ldots,a_n,b_1,\ldots,b_n),
    \end{equation}
    where
        \begin{align*}
        \calC_1 \coloneqq &  \biggl\{(a_1,\ldots,a_n)\in\Real^{n}:  a_1\geq \cdots\geq a_n\geq 0 \qq{and}
        \sum_{j=1}^n a_j = n\tau_i^2 \biggr\},
        \\
        \calC_2 \coloneqq & \biggl\{(b_1,\ldots,b_n)\in\Real^n: b_1\geq 0,\ldots,b_n\geq 0 \qq{and} \sum_{j=1}^n b_j \geq n \min_{\BFx'\in\calS}[k_i^0(\BFx,\BFx')]^2\biggr\}.
    \end{align*}
The reason for the inequality in \cref{eq:min-p-extend} is that the two minimization problems have the same objective function while the one in left-hand side has smaller feasible region.

    We now solve the minimization problem on the right-hand side of \cref{eq:min-p-extend}. Note that for any \((a_1,\ldots,a_n)\in\calC_1\), \(\min_{(b_1,\ldots,b_n)\in\calC_2} g(a_1,\ldots,a_n,b_1,\ldots,b_n)\) is a linear programming problem, and it is easy to see that its optimal solution is \(b_1^*=n \min_{\BFx'\in\calS} [k_i^0(\BFx,\BFx')]^2\) and \(b_2^*=\cdots=b_n^* =0\). Hence,
    \begin{equation}\label{eq:min-p-extend-solution}
        \min_{\substack{(a_1,\ldots,a_n)\in\calC_1 \\ (b_1,\ldots,b_n)\in \calC_2}}  g(a_1,\ldots,a_n,b_1,\ldots,b_n) =
        \min_{(a_1,\ldots,a_n)\in\calC_1} \frac{n \min_{\BFx'\in\calS} [k_i^0(\BFx,\BFx')]^2 }{a_1+\lambda_i^{\max}}
        =  \frac{n \min_{\BFx'\in\calS} [k_i^0(\BFx,\BFx')]^2 }{n \tau_i^2 +\lambda_i^{\max}}.
    \end{equation}
    Then, we can apply \cref{eq:min-p-extend,eq:min-p-extend-solution} to show \cref{eq:var_up_bound2}, completing the proof of \Cref{lemma:bound_acc_1}.
    \end{proof}

\begin{lemma}\label{lemma:bound_acc_2}
    Fix \(i=1,\ldots,M\) and a sampling policy \(\pi\). Suppose that the sequence of sampling locations \(\{\BFv^n:a^n=i,n=0,1,\ldots\}\) under \(\pi\) has an accumulation point \(\BFx_i^{\acc}\).
    If  \Cref{assump:stationary,assump:var,assump:compact} hold, then for any \(\epsilon>0\),
    \[\limsup_{n\to\infty}\max_{\BFx\in \calB(\BFx_i^{\acc}, \epsilon)}\Var^{\pi,n}[\theta_i(\BFx)] \leq \tau_i^2\qty[1-\rho_i^2(2\epsilon\BFone)],\]
    where \(\BFone\) is the vector of all ones with size $d \times 1$.
\end{lemma}
\begin{proof}[Proof of \Cref{lemma:bound_acc_2}.]
    It follows from \cref{eq:decreasing_var} that \(\{\Var^{\pi, n}[\theta_i(\BFx)]\}_{n=0}^\infty\) is a non-increasing sequence bounded below by zero. Hence, \(\Var^{\pi, n}[\theta_i(\BFx)]\) converges as \(n\to\infty\) and its limit is well defined.

    Fix \(\epsilon>0\). Let \(s_{i,n}\coloneqq \abs{\{\BFv^\ell\in \calB(\BFx_i^{\acc}, \epsilon):a^\ell=i,\ell=0,\ldots,n-1\}}\) be the number of times that alternative \(i\) is sampled at a point in \(\calB(\BFx_i^{\acc}, \epsilon)\) under \(\pi\) among the total \(n\) samples.
    Then, we must have \(s_{i,n}\to\infty\) since \(\BFx_i^{\acc}\) is an accumulation point. Note that reordering the sampling decision-observation pairs \(((a^0, \BFv^0),y^1), \ldots, ((a^{n-1}, \BFv^{n-1}),y^n)\) does not alter the conditional variance of \(\theta_i(\BFx)\). Hence, we may assume without loss of generality that the first \(s_{i,n}\) samples are all taken from alternative \(i\) at locations that belong to \(\calB(\BFx_i^{\acc},\epsilon)\). Since  the posterior variance decreases in the number of samples by \cref{eq:decreasing_var}, we conclude that for all \(\BFx \in \calB(\BFx_i^{\acc}, \epsilon)\),
    \vspace{-10pt} 
    \begin{equation}\label{eq:Scott}
        \Var^{\pi, n}[\theta_i(\BFx)] \leq \Var^{\pi, s_{i,n}}[\theta_i(\BFx)] \leq \tau_i^2 - \frac{ s_{i,n} \min_{\BFx'\in \calB(\BFx_i^{\acc}, \epsilon)} [k_i^0(\BFx,\BFx')]^2 }{s_{i,n} \tau_i^2 +\lambda_i^{\max}},
    \vspace{-5pt} 
    \end{equation}
    where the second inequality follows from  \Cref{lemma:bound_acc_1}.

    Note that \([k^0(\BFx,\BFx')]^2 = \tau_i^4 \rho_i^2(\abs{\BFx-\BFx'})\), and that  \(\norm{\BFx-\BFx'}\leq \norm{\BFx-\BFx_i^{\acc}}+ \norm{\BFx_i^{\acc}-\BFx'}\leq 2\epsilon\) for all \(\BFx,\BFx'\in \calB(\BFx_i^{\acc}, \epsilon)\).
    Hence, each component of \(\abs{\BFx-\BFx'} \) is no greater than \(2\epsilon\).
    Since \(\rho_i(\BFdelta)\) is decreasing in \(\BFdelta\) component-wise for \(\BFdelta>0\) (see \Cref{assump:stationary}), \(\rho_i(\abs{\BFx-\BFx'})\geq \rho_i(2\epsilon\BFone)\) for all \(\BFx,\BFx'\in \calB(\BFx_i^{\acc}, \epsilon)\).
    It then follows from  \cref{eq:Scott} that
    \vspace{-10pt} 
    \[\max_{\BFx\in \calB(\BFx_i^{\acc}, \epsilon)} \Var^{\pi, n}[\theta_i(\BFx)] \leq \tau_i^2 - \frac{s_{i,n} \tau_i^4\rho_i^2(2\epsilon\BFone)}{s_{i,n} \tau_i^2+\lambda_i^{\max}}. \]
    Sending \(n\to\infty\) completes the proof.
\end{proof}

\begin{lemma}\label{lemma:KG_bound}
Fix \(i=1,\ldots,M\).
If \Cref{assump:compact,assump:stationary,assump:var} hold, then for all \(\BFx\in\calX\),
    \[\IKG^n(i,\BFx) \leq \sqrt{\frac{ \tau_i^2\Var^n[\theta_i(\BFx)]}{2\pi \lambda_i(\BFx) c_i^2(\BFx)}}.\]
\end{lemma}
\begin{proof}[Proof of \Cref{lemma:KG_bound}.]
    Notice that
    \begin{align}
             & \int_{\calX}\E\qty[\max_{1\leq a\leq M} \mu^{n+1}_a(\BFv)\,\Big|\,\mathscr{F}^n, a^n=i, \BFv^n=\BFx]\gamma(\BFv) \dd{\BFv} \nonumber                                  \\
        = {}  & \int_{\calX}\E^n\qty[\max_{1\leq a\leq M} \qty(\mu^{n}_a(\BFv) + \sigma_a^n(\BFv,\BFv^n)Z^{n+1})\,\Big|\, a^n=i, \BFv^n=\BFx]\gamma(\BFv) \dd{\BFv}  \nonumber \\
        \leq {} & \int_{\calX}\max_{1\leq a\leq M} \mu^{n}_a(\BFv) \gamma(\BFv) \dd{\BFv} +
        \int_{\calX}\E^n\qty[\max_{1\leq a \leq M} \qty( \sigma_a^n(\BFv,\BFv^n)Z^{n+1})\,\Big|\,  a^n=i, \BFv^n=\BFx]\gamma(\BFv) \dd{\BFv}. \label{eq:upper_bound}
    \end{align}

    Since \(k_i^0(\BFx,\BFx')\) is a continuous function by \Cref{assump:stationary}, it follows from \Cref{assump:var} and the updating \cref{eq:mean_n} that \(\mu_i^n(\BFx)\) is a continuous function for any \(n\).
    Hence, \(\mu_i^n(\BFx)\) is bounded on \(\calX\) by \Cref{assump:compact}.
    This implies that the first integral in \cref{eq:upper_bound} is finite and can be subtracted from both sides of the inequality.
    Then, by the definition \cref{eq:int-KG},
    \begin{equation}\label{eq:IKG_upper_bound}
        \IKG^n(i,\BFx) \leq \frac{1}{c_i(\BFx)} \int_{\calX}\E^n\qty[\max_{1\leq a \leq M} \qty( \sigma_a^n(\BFv,\BFv^n)Z^{n+1})\,\Big|\,  a^n=i, \BFv^n=\BFx]\gamma(\BFv) \dd{\BFv} \coloneqq I.
    \end{equation}
    It follows from \cref{eq:sigma} that
    \begin{align}
        I & = \frac{1}{c_i(\BFx)} \int_{\calX}\E^n\qty[\max\qty{\sigma_{a^n}^n(\BFv,\BFv^n)Z^{n+1}, 0 }  \,\Big|\,  a^n=i, \BFv^n=\BFx]\gamma(\BFv) \dd{\BFv} \nonumber \\[0.5ex]
          & = \frac{1}{c_i(\BFx)} \int_{\calX}\E^n\qty[\max\qty{\abs{\tilde\sigma_{i}^n(\BFv,\BFx)}Z^{n+1}, 0 } ]\gamma(\BFv) \dd{\BFv} \nonumber                           \\[0.5ex]
          & = \frac{1}{c_i(\BFx)} \int_{\calX}\E^n \qty[\abs{\tilde\sigma_i^n(\BFv,\BFx)}Z^{n+1}\ind_{\{\abs{\tilde\sigma_i^n(\BFv,\BFx)} Z^{n+1}> 0\}}]\gamma(\BFv) \dd{\BFv} \nonumber \\[0.5ex]
          & = \frac{1}{c_i(\BFx)} \int_{\calX}\abs{\tilde\sigma_i^n(\BFv,\BFx)} \qty[ \int_0^\infty z \phi(z) \dd{z}]\gamma(\BFv) \dd{\BFv} \nonumber                                 \\[0.5ex]
          & = \frac{1}{\sqrt{2\pi} c_i(\BFx)} \int_{\calX}\abs{\tilde\sigma_i^n(\BFv,\BFx)}\gamma(\BFv) \dd\BFv. \label{eq:calculating_I}
    \end{align}
    Moreover, by \cref{eq:sigma},
    \begin{equation}\label{eq:sigma_upper_bound}
        \abs{\tilde\sigma_i^n(\BFv,\BFx)} = \abs{\frac{\Cov^n[\theta_i(\BFv),\theta_i(\BFx)]}{\sqrt{\Var^n[\theta_i(\BFx)]+\lambda_i(\BFx)}}} \leq \sqrt{\frac{\Var^n[\theta_i(\BFv)]\Var^n[\theta_i(\BFx)]}{\Var^n[\theta_i(\BFx)]+\lambda_i(\BFx)}} \leq \sqrt{\frac{\tau_i^2\Var^n[\theta_i(\BFx)]}{\lambda_i(\BFx)}},
    \end{equation}
    where the last inequality follows because \(0 \leq \Var^n[\theta_i(\BFv)]\leq \Var[\theta_i(\BFv)]=\tau_i^2\) for all \(\BFv\in\calX\) by \cref{eq:decreasing_var}.
    The proof is completed by combining \cref{eq:IKG_upper_bound,eq:sigma_upper_bound,eq:calculating_I}.
\end{proof}

\begin{lemma}\label{lemma:liminf}
    Fix \(i=1,\ldots,M\). If \Cref{assump:stationary,assump:compact,assump:var} hold, and \(\eta_i^\infty=\infty\) under the IKG policy, then for any \(\BFx\in\calX\),
    \[\liminf_{n\to\infty} \IKG^n(i,\BFx) = 0.\]
\end{lemma}
\begin{proof}[Proof of \Cref{lemma:liminf}.]
    Since \(\calX\) is compact by \Cref{assump:compact}, the sequence \(\{\BFv^n\in\calX:a^n=i,n=0,1,\ldots,\}\) is bounded, and it is of length \(\eta_i^\infty=\infty\).
    Hence, it has an accumulation point \(\BFx_i^{\acc}\).
    Let \(\{\ell_{i,n}\}_{n=0}^\infty\) be the subsequence of \(\{n:a^n=i, n=0,1,\ldots\}\) such that \(\ell_{i,n}\to\infty\) and \(\BFv^{\ell_{i,n}}\to \BFx_i^{\acc}\) as \(n\to\infty\).
    Fix \(\epsilon>0\).
    Then, by \Cref{lemma:bound_acc_2},
    \begin{equation*}
        \limsup_{n\to\infty}\Var^n[\theta_i(\BFv^{\ell_{i,n}})] \leq \tau_i^2[1-\rho_i^2(2\epsilon\BFone)].
    \end{equation*}
    It then follows from \Cref{lemma:KG_bound} that
    \begin{equation*}
        \limsup_{n\to\infty} \IKG^{\ell_{i,n}}(i,\BFv^{\ell_{i,n}})  \leq
        \limsup_{n\to\infty} \sqrt{\frac{\tau_i^2\Var^{\ell_{i,n}}[\theta_i(\BFv^{\ell_{i,n}})]}{2\pi \lambda_i(\BFx) c_i^2(\BFx)}}
        \leq \sqrt{\frac{\tau_i^4[1-\rho_i^2(2\epsilon\BFone)]}{2\pi \lambda_i(\BFx) c_i^2(\BFx)}}.
    \end{equation*}
    By sending \(\epsilon\to 0\), we have \(\rho_i(2\epsilon\BFone)\to 1\) and thus,
    \(\limsup_{n\to\infty} \IKG^{\ell_{i,n}}(i,\BFv^{\ell_{i,n}})  \leq  0\).
    Since the limit inferior of a sequence is no greater than that of its subsequence,
    \begin{equation}\label{eq:liminf_negative}
        \liminf_{n\to\infty} \IKG^{n}(i,\BFv^n)  \leq  \liminf_{n\to\infty} \IKG^{\ell_{i,n}}(i,\BFv^{\ell_{i,n}}) \leq  \limsup_{n\to\infty} \IKG^{\ell_{i,n}}(i,\BFv^{\ell_{i,n}}) \leq 0.
    \end{equation}

    Moreover, by the definition of IKG \cref{eq:int-KG} and Jensen's inequality,
    \[
        \IKG^n(i,\BFx)  \geq \frac{1}{c_i(\BFx)} \int_{\calX}\qty{\max_{1\leq a\leq M}\E\qty[ \mu^{n+1}_a(\BFv)\,\Big|\,\mathscr{F}^n, a^n=i, \BFv^n=\BFx] - \max_{1\leq a\leq M} \mu^n_a(\BFv)} \gamma(\BFv) \dd{\BFv} = 0,
    \]
    for each \(i=1,\ldots,M\) and \(\BFx\in\calX\), where the equality follows immediately from the updating \cref{eq:mean_update}.
    This, in conjunction with \cref{eq:liminf_negative}, implies that
    \(\liminf_{n\to\infty} \IKG^{n}(i,\BFv^n)  = 0\).
    By the definition of the sampling location \(\BFv^n\) in \cref{eq:max-int-KG},
    \(\IKG^{n}(i,\BFv^n) = \max_{\BFx\in\calX}\IKG^{n}(i,\BFx)\).
    Hence, for any \(\BFx\in\calX\),
\begin{equation}\label{eq:liminfIKG}
0 \leq \liminf_{n\to\infty} \IKG^n(i,\BFx) \leq\liminf_{n\to\infty} \IKG^{n}(i,\BFv^n) = 0,
\end{equation}
    which completes the proof.
\end{proof}

Before we prove \Cref{prop:step1}, we need one more technical result about the almost sure convergence of $\mu_i^n$, which is stated in the following \Cref{lemma:u_n_converge}.
A similar result is also given in \citet[Proposition 2.9]{bect2019}, and the proof there directly applies here.

\begin{lemma}\label{lemma:u_n_converge}
If \Cref{assump:stationary,assump:compact,assump:var} hold, then for all \(i=1,\ldots,M\), \(\mu_i^n(\BFx)\) converges to $\mu_i^\infty(\BFx)\coloneqq \E[\theta_i(\BFx)|\mathscr{F}^\infty]$ uniformly in \(\BFx\in\calX\) a.s.  as \(n\to\infty\). That is,
\[
\pr\left\{\omega: \sup_{\BFx\in\calX} |\mu_i^n(\BFx;\omega) - \mu_i^\infty(\BFx;\omega)|\to 0 \right\}=1.\]
\end{lemma}

\begin{proof}[Proof of \Cref{lemma:u_n_converge}.]
Fix \(i=1,\ldots,M\).
Note that $\theta_i$ is a Gaussian process under the prior.
It follows from \Cref{assump:stationary,assump:var} and Theorem 1.4.1 of~\cite{AdlerTaylor07} that the sample paths of $\theta_i$ are continuous a.s..
The proof is then completed by directly following the arguments in the proof of Proposition 2.9 in \cite{bect2019}.
\end{proof}

We are now ready to prove \Cref{prop:step1}.
\begin{proof}[Proof of \Cref{prop:step1}.]

Let \(\mu^n \coloneqq (\mu_1^n,\ldots,\mu_M^n)\) denote the posterior mean of \((\theta_1,\ldots,\theta_M)\) conditioned on \(\mathscr{F}^n\).
Let \(\omega\) denote a generic sample path.
Fix \(i=1,\ldots,M\).
Define
$$
\Omega_0 = \{\omega: \eta_i^\infty(\omega) = \infty, \ \mu^n(\BFx;\omega) \rightarrow \mu^\infty(\BFx; \omega) \text{ uniformly in $\BFx\in\calX$ as }n\to\infty\}.
$$
Then, \(\pr(\Omega_0)=1\) by the assumption of \Cref{prop:step1} and \Cref{lemma:u_n_converge}.
Fix an arbitrary \(\BFx\in\calX\).
We now prove that, under the IKG policy,
\begin{equation}\label{eq:omega0_result}
    k_i^\infty(\BFx,\BFx; \omega)=0, \text{ for any }\omega \in \Omega_0,
\end{equation}
which establishes \Cref{prop:step1}.
We prove \cref{eq:omega0_result} by contraction and assume that there exists some \(\omega_0 \in \Omega_0\) such that \(k_i^\infty(\BFx,\BFx;\omega)>0\).
In the remaining proof, we suppress the sample path \(\omega_0\) to simplify notation.

It follows from the continuity of \(k_i^0(\BFx,\cdot)\) assumed in \Cref{assump:stationary} and the updating \cref{eq:cov_n} that \(k_i^n(\BFx,\cdot)\) is continuous.
The uniform convergence of \(k_i^n(\BFx,\cdot)\) by \Cref{prop:unif_conv} then implies that \(k^\infty_i(\BFx,\cdot)\) is also continuous.
Hence, there exist \(\epsilon>0\) such that \(\min_{\BFv\in\calB(\BFx, \epsilon)}k_i^\infty(\BFx,\BFv)>0\).
The uniform convergence of \(k^n_i(\BFx,\cdot)\) further implies that there exists \(\delta>0\) such that \(k^n_i(\BFx,\BFv)\geq \delta\) for all \(\BFv\in\calB(\BFx, \epsilon)\) and \(n\geq 1\).
By \cref{eq:sigma},
\begin{align*}\label{eq:unif_lower_bound}
    \inf_{\BFv\in\calB(\BFx, \epsilon),n\geq 1}\tilde\sigma_i^n(\BFv,\BFx)
     & = [k_i^n(\BFx,\BFx)+\lambda_i(\BFx)]^{-1/2}\inf_{\BFv\in\calB(\BFx, \epsilon),n\geq 1}k_i^n(\BFv,\BFx) \\
     &\geq \delta[k_i^0(\BFx,\BFx)+\lambda_i(\BFx)]^{-1/2} \coloneqq \alpha_1>0.
\end{align*}

Let \(g(s, t)\coloneqq t\phi(s/t)-s\Phi(-s/t)\); see \Cref{lemma:func_g} for properties of \(g(s,t)\), including positivity and monotonicity.
Then,
\[\abs{\tilde\sigma_i^n(\BFv, \BFx)}\phi\qty(\abs{\frac{\Delta_i^n(\BFv)}{\tilde\sigma_i^n(\BFv, \BFx)}})
    -\abs{\Delta_i^n(\BFv)}\Phi\qty(-\abs{\frac{\Delta_i^n(\BFv)}{\tilde\sigma_i^n(\BFv, \BFx)}})=g(\abs{\Delta_i^n(\BFv)}, \abs{\tilde\sigma_i^n(\BFv, \BFx)}) \geq 0,
\]
for all \(\BFv\in\calX\).
Consequently, \Cref{lemma:calculating_IKG} implies that
\begin{equation*}
    \IKG^n(i,\BFx) \geq \frac{1}{c_i(\BFx)} \int_{\calB(\BFx, \epsilon)}g(\abs{\Delta_i^n(\BFv)}, \abs{\tilde\sigma_i^n(\BFv, \BFx)})\gamma(\BFv) \dd{\BFv}
    \geq \frac{1}{c_i(\BFx)} \int_{\calB(\BFx, \epsilon)} g(\abs{\Delta_i^n(\BFv)},\alpha_1)\gamma(\BFv) \dd\BFv,
\end{equation*}
for all \(n\geq 1\), where the second inequality holds because \(g(s,t)\) is strictly increasing in \(t\in(0,\infty)\).
Note that \(\liminf_{n\to\infty} \IKG^n(i,\BFx)=0\)  by \Cref{lemma:liminf}.
Hence,
\begin{equation}\label{eq:IKG_lower_bound}
    0 \geq \liminf_{n\to\infty} \frac{1}{c_i(\BFx)} \int_{\calB(\BFx, \epsilon)} g(\abs{\Delta_i^n(\BFv)},\alpha_1) \gamma(\BFv) \dd\BFv \geq \frac{1}{c_i(\BFx)} \int_{\calB(\BFx, \epsilon)} \liminf_{n\to\infty} g(\abs{\Delta_i^n(\BFv)},\alpha_1) \gamma(\BFv) \dd\BFv,
\end{equation}
where the second inequality holds due to Fatou's lemma.
Furthermore, since for any \(\BFv\in\calX\), \(\Delta_i^n(\BFv)=\mu_i^n(\BFv)-\max_{a\neq i}\mu_a^n(\BFv)\), and \(\mu_a^n(\BFv) \rightarrow \mu_a^\infty(\BFv)\) for \(a=1,\ldots,M\), where \(\abs{\mu_a^\infty(\BFv)} < \infty\).
then
\[\limsup_{n\to\infty}\abs{\Delta_i^n(\BFv)}  \leq \limsup_{n\to\infty} [2\max_a\abs{\mu^n_a(\BFv)}]=  2\max_a\abs{\mu^\infty_a(\BFv)}\coloneqq \alpha_2(\BFv)<\infty,\]
for all \(\BFv\in\calB(\BFx, \epsilon)\).
Then, in the light of \cref{eq:IKG_lower_bound} and the fact that \(g(s,t)\) is strictly decreasing in \(s\in[0,\infty)\),
\[0\geq \frac{1}{c_i(\BFx)} \int_{\calB(\BFx, \epsilon)}   \liminf_{n\to\infty} g(\abs{\Delta_i^n(\BFv)},\alpha_1) \gamma(\BFv) \dd{\BFv}
   \geq \frac{1}{c_i(\BFx)} \int_{\calB(\BFx, \epsilon)}  g(\alpha_2(\BFv), \alpha_1) \gamma(\BFv) \dd{\BFv}.\]
This contracts the fact that \(g(s,t)>0\) for all \(s\in[0,\infty)\) and \(t\in(0,\infty)\).
Therefore, \cref{eq:omega0_result} is proved.
\end{proof}

\subsection{D. Proof of Proposition \ref{prop:step2}}\label{sec:step2}

Let \(S^n\coloneqq (\mu_1^n,\ldots,\mu_M^n,k_1^n,\ldots,k_M^n)\) denote the state at time \(n\), which fully determines the posterior distribution of \((\theta_1,\ldots,\theta_M)\) conditioned on \(\mathscr{F}^n\).
The state transition \(S^n \to S^{n+1}\) is governed by \cref{eq:mean_update,eq:cov_update}, which is determined by the sampling decision \((a^n,\BFv^n)\).

Let \(s\coloneqq(\mu_1,\ldots,\mu_M,k_1,\ldots,k_M)\in\mathbb{S}\) be a generic state and \(\mathbb{S}\) denote the set of states for which \(\mu_i\) is a continuous function and \(k_i\) is a continuous covariance function for each \(i=1,\ldots,M\).
For \(s\in\mathbb{S}\), define
\[V(s)\coloneqq \int_{\calX}\max_{1\leq a\leq M} \mu_a(\BFv) \gamma(\BFv) \dd\BFv,\]
and
\[Q(s,i,\BFx)\coloneqq \E[V(S^{n+1})\,|\,S^n=s, a^n=i,\BFv^n=\BFx].\]
By the following \Cref{lemma:Qfactors-IKG}, it is easy to see that at time \(n\), the IKG policy \eqref{eq:max-int-KG} chooses
\begin{equation}\label{eq:IKG-3}
    \argmax_{1\leq i\leq M,\BFx\in\calX} \ [c_i(\BFx)]^{-1} [Q(S^n,i,\BFx) - V(S^n)].
\end{equation}

\begin{lemma}\label{lemma:Qfactors-IKG}
    Fix \(s\in\mathbb{S}\), \(i=1,\ldots,M\), and \(\BFx\in\calX\) where \(\calX\) is compact,
    \[
        Q(s,i,\BFx) = \int_{\calX}\E\qty[\max_{1\leq a\leq  M}\mu_a^{n+1}(\BFv)\,\Big|\, S^n=s, a^n=i,\BFv^n=\BFx] \gamma(\BFv) \dd\BFv.
    \]
\end{lemma}
\begin{proof}[Proof of \Cref{lemma:Qfactors-IKG}]
    Notice that by the updating \cref{eq:mean_update}, given \(S^n=s\), \(a^n=i\), \(\BFv^n=\BFx\) and \(Z^{n+1}\),
    \[
        \max_{1\leq a\leq  M}\mu_a^{n+1}(\BFv) = \max\qty{\mu_i(\BFv)+ \tilde\sigma_i(\BFv,\BFx)Z^{n+1},\max_{a\neq i}\mu_a(\BFv)}.
    \]
    Let \(f\qty(s,i,\BFx,\BFv,Z^{n+1}) = \max_{1\leq a\leq  M}\mu_a^{n+1}(\BFv) - \max_{a\neq i}\mu_a(\BFv)\).
    Then \(f\qty(s,i,\BFx,\BFv,Z^{n+1}) \geq 0\), for all \(\BFv \in \calX\) and \(Z^{n+1}\).
    Hence,
    \begin{align*}
        Q(s,i,\BFx) & = \E\qty[\int_{\calX}\max_{1\leq a\leq  M}\mu_a^{n+1}(\BFv)\gamma(\BFv)\dd\BFv\,\Big|\, S^n=s, a^n=i,\BFv^n=\BFx] \nonumber                   \\[0.5ex]
                    & = \E\qty[\int_{\calX} \Big(f\qty(s,i,\BFx,\BFv,Z^{n+1}) + \max_{a\neq i}\mu_a(\BFv) \Big) \gamma(\BFv) \dd{\BFv}]                          \\[0.5ex]
                    & = \E\qty[\int_{\calX} f\qty(s,i,\BFx,\BFv,Z^{n+1}) \gamma(\BFv) \dd{\BFv} ] + \int_{\calX} \max_{a\neq i}\mu_a(\BFv) \gamma(\BFv) \dd{\BFv} \\[0.5ex]
                    & =\int_{\calX} \E\qty[f\qty(s,i,\BFx,\BFv,Z^{n+1}) ]\gamma(\BFv) \dd{\BFv} + \int_{\calX} \max_{a\neq i}\mu_a(\BFv) \gamma(\BFv) \dd{\BFv},
    \end{align*}
    where the interchange of integral and expectation is justified by Tonelli's theorem for nonnegative functions, and \(\int_{\calX} \max_{a\neq i}\mu_a(\BFv) \gamma(\BFv) \dd{\BFv}\) is finite since \(\mu_i(\BFv)\) is continuous on the compact set \(\calX\) for \(i=1,\ldots,M\).
    Thus the result in \Cref{lemma:Qfactors-IKG} follows immediately.
\end{proof}

\begin{lemma}\label{lemma:Qfactors}
    Fix \(s\in\mathbb{S}\), \(i=1,\ldots,M\), and \(\BFx\in\calX\) where \(\calX\) is compact and \(\alpha(\cdot) >0\) on \(\calX\).
    Then, \(Q(s,i,\BFx)\geq V(s)\) and the equality holds if and only if \(k_i(\BFx,\BFx)= 0\).
\end{lemma}
\begin{proof}[Proof of \Cref{lemma:Qfactors}.]
    Applying \Cref{lemma:Qfactors-IKG} and the updating \cref{eq:mean_update},
    \begin{align}
        Q(s,i,\BFx)
         & = \int_{\calX}\E\qty[\max\qty{\mu_i(\BFv)+ \tilde\sigma_i(\BFv,\BFx)Z^{n+1},\max_{a\neq i}\mu_a(\BFv)}]\gamma(\BFv)\dd{\BFv}\nonumber             \\[0.5ex]
         & \geq \int_{\calX} \max\qty{\E\qty[\mu_i(\BFv)+ \tilde\sigma_i(\BFv,\BFx)Z^{n+1}],\max_{a\neq i}\mu_a(\BFv)}\gamma(\BFv)\dd\BFv \label{eq:jensen} \\[0.5ex]
         & = \int_{\calX}\max_{1\leq a\leq M}  \mu_a(\BFv)\gamma(\BFv)\dd\BFv = V(s),\nonumber
    \end{align}
    where \cref{eq:jensen} follows from Jensen's inequality since \(\max(\cdot,\cdot)\) is a strictly convex function.

    If \(k_i(\BFx,\BFx)=0\), then in the light of the fact that \(k_i\) is a covariance function, we must have that
    \[\abs{k_i(\BFv,\BFx)} \leq \sqrt{k_i(\BFv,\BFv)k_i(\BFx,\BFx)} = 0,\]
    so \(k_i(\BFv,\BFx)=0\) for all \(\BFv\in\calX\). Hence,  \(\tilde\sigma_i^n(\BFv,\BFx)= 0\) by \cref{eq:sigma}, so
    \(\mu_a^{n+1}(\BFv) = \mu_a^n(\BFv)\) for all \(a=1,\ldots,M\) and \(\BFv\in\calX\). Hence, \(\mu_a^{n+1}(\BFv)\) is deterministic given \(S^n\)  for all \(a=1,\ldots,M\) and \(\BFv\in\calX\). Thus, the inequality \cref{eq:jensen} holds with equality.

    Next, assume conversely that \(Q(s,i,\BFx)=V(s)\). If \(k_i(\BFx,\BFx)\neq 0\), then the continuity of \(k_i\) implies that \(k_i(\BFv,\BFx)\neq 0\) for all \(\BFv\in\tilde\calX\), where \(\tilde\calX\subset\calX\) is an open neighborhood of \(\BFx\).
    Without loss of generality, we assume that for all \(\BFv\in\tilde\calX\), \(k_i(\BFv,\BFx)>0\) and thus, \(\tilde\sigma_i(\BFv,\BFx)=k_i(\BFv,\BFx)/\sqrt{k_i(\BFx,\BFx)+\lambda_i(\BFx)}>0\). By the strict convexity of \(\max(\cdot,\cdot)\) and Jensen's inequality,
    \[\E\qty[\max\qty{\mu_i(\BFv)+ \tilde\sigma_i(\BFv,\BFx)Z^{n+1},\max_{a\neq i}\mu_a(\BFv)}]
        > \max\qty{\E\qty[\mu_i(\BFv)+ \tilde\sigma_i(\BFv,\BFx)Z^{n+1}],\max_{a\neq i}\mu_a(\BFv)},\]
    for \(\BFv\in\tilde\calX\). Hence,  \cref{eq:jensen} becomes a strict inequality since \(\gamma(\BFv)>0\) for all \(\BFv\in\calX\).
    This contradicts \(Q(s,i,\BFx)=V(s)\), so \(k_i(\BFx,\BFx)= 0\).
\end{proof}

\begin{lemma}\label{lemma:limiting_cov}
    Fix \(i=1,\ldots,M\). If \Cref{assump:stationary,assump:var} hold, and \(k_i^\infty(\BFx,\BFx)=0\) for some \(\BFx\in\calX\), then \(\eta_i^\infty=\infty\).
\end{lemma}
\begin{proof}[Proof of \Cref{lemma:limiting_cov}.]
    We prove by contradiction and assume that \(\eta_i^\infty<\infty\).
    Then, \(N_i\coloneqq \min\{n:\eta_i^n=\eta_i^\infty\} <\infty\) and \(a^n\neq i\) for all \(n\geq N_i\).
    Due to the mutual independence between the alternatives, it follows that the posterior distribution of \(\theta_i\) remains the same for \(n\geq N_i\).
    In particular, \(k_i^n(\BFx,\BFx)=k_i^{N_i}(\BFx,\BFx)\) for all \(n>N_i\). Hence, \(k_i^{N_i}(\BFx,\BFx)=k_i^\infty(\BFx,\BFx)=0\).
    It follows from \cref{eq:cov_update} that
    \[k_i^{N_i-1}(\BFx,\BFx) = k_i^{N_i}(\BFx,\BFx) + \qty[\sigma_i^{N_i-1}(\BFx,\BFv^{N_i-1})]^2
        = \qty[\sigma_i^{N_i-1}(\BFx,\BFv^{N_i-1})]^2.
    \]
    By the definition of \(N_i\), \(a^{N_i-1} = i\).
    Then by \cref{eq:sigma},
    \begin{equation}\label{eq:recursion_post_var}
        k_i^{N_i-1}(\BFx,\BFx) = \frac{\qty[k_i^{N_i-1}(\BFx,\BFv^{N_i-1})]^2}{k_i^{N_i-1}(\BFv^{N_i-1},\BFv^{N_i-1})+\lambda_i(\BFv^{N_i-1})}.
    \end{equation}
    Notice that
    \begin{align}
        \nonumber                                                                                                                      \\[-7ex]
        [k_i^{N_i-1}(\BFx,\BFv^{N_i-1})]^2 & = \qty{\Cov^{N_i-1}\qty[\theta_i(\BFx),\theta_i(\BFv^{N_i-1})]}^2 \nonumber  \\[0.5ex]
                                           & \leq \Var^{N_i-1}[\theta_i(\BFx)]  \Var^{N_i-1}[\theta_i(\BFv^{N_i-1})] \nonumber         \\[0.5ex]
                                           & = k_i^{N_i-1}(\BFx,\BFx)k_i^{N_i-1}(\BFv^{N_i-1},\BFv^{N_i-1}). \label{eq:Cauchy_Schwarz}
    \end{align}
    It follows from \cref{eq:recursion_post_var,eq:Cauchy_Schwarz} that \(\lambda_i(\BFv^{N_i-1}) k_i^{N_i-1}(\BFx,\BFx)\leq 0\).
Thus, \(k_i^{N_i-1}(\BFx,\BFx)=0\), since \(k_i^{N_i-1}(\BFx,\BFx)\geq 0\) and \(\lambda_i(\BFv^{N_i-1})>0\) in \Cref{assump:var}.
    By induction, we can conclude that \(k_i^0(\BFx,\BFx)=0\), which contracts the fact that \(k_i^0(\BFx,\BFx)=\tau_i^2>0\) in \Cref{assump:stationary}.
    Therefore, we must have \(\eta_i^\infty=\infty\).
\end{proof}

We are now ready to prove \Cref{prop:step2}.
\begin{proof}[Proof of \Cref{prop:step2}.]
    Define \(\Omega_1\coloneqq \{\omega: S^n(\omega) \to S^\infty(\omega)\text{ pointwise as }n\to\infty\}\).
    By  \Cref{lemma:u_n_converge} and \Cref{prop:unif_conv}, \(\pr(\Omega_1)=1\).
    For any \(i=1,\ldots,M\), define the event \(\SFH_i \coloneqq\{\omega: \eta_i^\infty(\omega)<\infty\}\). Then, \(k_i^\infty(\BFx,\BFx;\omega)>0\) for all \(\BFx\in\calX\) and \(\omega\in\SFH_i\) by \Cref{lemma:limiting_cov}.
    On the other hand, \Cref{prop:step1} implies that \(k_i^\infty(\BFx,\BFx;\omega)=0\) for all \(\BFx\in\calX\) and \(\omega\in\SFH_i^{\mathsf{c}}\cap\Omega_1\), where \(\SFH_i^{\mathsf{c}}\) is the complement of \(\SFH_i\).
    Thus, by \Cref{lemma:Qfactors},
    \begin{equation}\label{eq:yesno_benefit}
        \begin{aligned}
             & Q(S^\infty(\omega), i,\BFx) > V(S^\infty(\omega)),\quad \text{for all }\omega \in \SFH_i \cap\Omega_1,             \\
             & Q(S^\infty(\omega), i,\BFx) = V(S^\infty(\omega)),\quad \text{for all }\omega \in \SFH_i^{\mathsf{c}}\cap\Omega_1.
        \end{aligned}
    \end{equation}
    Further, for any subset \(A\subseteq \{1,\ldots,M\}\), define the event
    \[
        \SFH_A\coloneqq\{\cap_{i\in A}\SFH_i\} \cap \{\cap_{i\notin A}\SFH_i^{\mathsf{c}}\}.
    \]

    Choose any \(A\neq \emptyset\).
    When \(A = \{1,\ldots,M\}\), \(\SFH_A=\emptyset\), because it is impossible that all alternative have finite samples while \(n\to\infty\).
    So \(\SFH_A\cap \Omega_1 = \emptyset\).
    When \(A \neq \{1,\ldots,M\}\), we prove \(\SFH_A\cap \Omega_1 = \emptyset\) by contradiction.
    Suppose that  \( \SFH_A\cap \Omega_1 \neq \emptyset\) so that we can choose and fix a sample path \(\omega_0 \in \SFH_A\cap \Omega_1\).
    Then, \(\eta_i^\infty(\omega_0)<\infty\) for all \(i \in A\).
    Hence, there exists \(T_i(\omega_0)<\infty\) for all \(i\in A\) such that the IKG policy does not choose alternative \(i\) for \(n > T_i(\omega_0)\).
    Let \(T(\omega_0)\coloneqq \max_{i\in A}T_i(\omega_0)\).
    Then, \(T(\omega_0)<\infty\) and the IKG policy does not choose \(i\in A\) for \(n > T(\omega_0)\).
    On the other hand, it follows from \cref{eq:yesno_benefit} that for all \(i\in A\), \(i'\notin A\), and \(\BFx\in\calX\),
    \[Q(S^\infty(\omega_0), i,\BFx) - V(S^\infty(\omega_0)) > Q(S^\infty(\omega_0), i',\BFx) - V(S^\infty(\omega_0)) = 0.\]
    Let \(Q^\dag(s, i,\BFx) \coloneqq Q(s, i,\BFx) - V(s)\) for simplicity.
    Then, by virtue of the compactness of \(\calX\) and the positivity of \(c_i(\BFx)\),
    \[\max_{\BFx\in\calX} \ [c_i(\BFx)]^{-1} Q^\dag(S^\infty(\omega_0), i,\BFx)
        > \max_{\BFx\in\calX} \ [c_{i'}(\BFx)]^{-1} Q^\dag(S^\infty(\omega_0), i',\BFx)=0,\]
    for all \(i\in A\) and \(i'\notin A\).
    Hence,
    \begin{equation}\label{eq:contradiction_A}
    \min_{i\in A} \max_{\BFx\in\calX} \ [c_i(\BFx)]^{-1} Q^\dag(S^\infty(\omega_0), i,\BFx)
        > \max_{i'\notin A} \max_{\BFx\in\calX} \ [c_{i'}(\BFx)]^{-1} Q^\dag(S^\infty(\omega_0), i',\BFx)=0.
    \end{equation}
    Notice that \(S^n(\omega_0) \to S^\infty(\omega_0)\) pointwise as \(n\to\infty\) since \(\omega_0\in\Omega_1\).
    Hence, there exists a finite number \(\tilde n(\omega_0) > T(\omega_0)\) such that
    \[\min_{i\in A}\max_{\BFx\in\calX} \ [c_i(\BFx)]^{-1} Q^\dag(S^{\tilde n(\omega_0) }(\omega_0), i,\BFx)
        > \max_{i'\notin A} \max_{\BFx\in\calX} \ [c_{i'}(\BFx)]^{-1} Q^\dag(S^{\tilde n(\omega_0) }(\omega_0), i',\BFx),\]
    which implies that IKG policy must choose alternative \(i\in A\) at time \(\tilde n(\omega_0)\) by \cref{eq:IKG-3}.
    This contradicts the definition of \( T(\omega_0)\).
    Therefore, the event \( \SFH_A\cap \Omega_1\) must be empty for any nonempty \(A\subseteq \{1,\ldots,M\}\).

    It then follows immediately that \(\pr(\SFH_A)=0\) for any nonempty \(A\subseteq \{1,\ldots,M\}\), since \(\pr(\Omega_1)=1\).
    Notice that the whole sample space \(\Omega = \cup_{A\subseteq \{1,\ldots,M\}} \SFH_A\).
    Hence,
    \[
        1 =\pr(\SFH_{\emptyset}) = \pr\qty(\cap_{i = 1}^M \SFH_i^{\mathsf{c}} )
        = \pr \qty( \{\omega:\eta_i^\infty=\infty \text{ for all } i=1,\ldots,M \} ),
    \]
    which completes the proof.
\end{proof}

\subsection{E. Proof of Theorem \ref{theo:consistency}}

\begin{proof}[\unskip\nopunct]
    Part (i) is an  immediate consequence of  \Cref{prop:step1,prop:step2}.
    The other two parts follow closely the proof of similar results in Theorem 1 of~\cite{XieFrazierChick16}.

    For part (ii), fix an arbitrary \(\BFx\in\calX\). Note that for each \(i=1,\ldots,M\),
    \[\E[(\mu^n_i(\BFx)-\theta_i(\BFx))^2]=\E[\E^n[(\mu^n_i(\BFx)-\theta_i(\BFx))^2]]=\E[k^n_i(\BFx,\BFx)]\to \E[k_i^\infty(\BFx,\BFx)]=0,\]
    as \(n\to\infty\), where the convergence holds due to the fact that \(0\leq k^n_i(\BFx,\BFx)\leq k^0_i(\BFx,\BFx)\) from \cref{eq:decreasing_var} and the dominated convergence theorem.
    This asserts that \(\mu^n_i(\BFx)\to\theta_i(\BFx)\) in \(\mathsf{L}^2\).
    By \Cref{lemma:u_n_converge}, \(\mu_i^n(\BFx)\to\mu_i^\infty(\BFx)\) a.s., which implies that \(\theta_i(\BFx)=\mu_i^\infty(\BFx)\) a.s., due to the a.s. uniqueness of convergence in probability.
    Thus, \(\mu^n_i(\BFx)\to\theta_i(\BFx)\) a.s. as \(n\to\infty\).

    For part (iii), let us again fix \(\BFx\in\calX\). Let \(i^*(\BFx)\in\argmax_i\theta_i(\BFx)\).
    We now show that \(\argmax_i\mu_i^n(\BFx)\to i^*(\BFx)\) a.s. as \(n\to\infty\).
    Again, we let \(\omega\) denote a generic sample path and use notations like \(i^*(\BFx;\omega)\) to emphasize the dependence on \(\omega\).
    Let \(\epsilon(\BFx;\omega)\coloneqq \theta_{i^*(\BFx,\omega)}(\BFx;\omega) - \max_{a\neq i^*(\BFx,\omega)}\theta_a(\BFx;\omega)\).
    Then, \(\pr(\{\omega:\epsilon(\BFx;\omega)>0\})=1\) because \((\theta_1(\BFx;\omega),\ldots,\theta_M(\BFx;\omega))\) is a realization of a multivariate normal random variable under the prior distribution.
    Hence, the event \(\tilde\Omega\coloneqq \{\omega:\epsilon(\BFx;\omega)>0 \text{ and } \mu^n_i(\BFx;\omega)\to\theta_i(\BFx;\omega) \text{ for all }\) \(i = 1,\ldots,M\}\) occurs with probability 1.
    Fix an arbitrary \(\tilde\omega\in\tilde\Omega\).
    To complete the proof, it suffices to show that \(\argmax_i \mu^n_i(\BFx;\tilde\omega)\to i^*(\BFx;\tilde\omega)\) as \(n\to\infty\).

    Clearly, there exists \(N(\tilde\omega)<\infty\) such that
    \(\abs{\mu_i^n(\BFx;\tilde\omega)-\theta_i(\BFx;\tilde\omega)}<\epsilon(\BFx;\tilde\omega)/2\) for all \(n>N(\tilde\omega)\) and \(i=1,\ldots,M\).
    Hence, for all \(i\neq i^*(\BFx;\tilde\omega)\) and \(n\geq N(\tilde\omega)\),
    \[\mu^n_{i^*(\BFx;\tilde\omega)}(\BFx;\tilde\omega) >\theta_{i^*(\BFx;\tilde\omega)}(\BFx;\tilde\omega) - \frac{\epsilon(\BFx;\tilde\omega)}{2}
        \geq \theta_i(\BFx;\tilde\omega)+\frac{\epsilon(\BFx;\tilde\omega)}{2}>\mu_i^n(\BFx;\tilde\omega).\]
    This implies that \(i^*(\BFx;\tilde\omega) = \argmax_i \mu_i^n(\BFx;\tilde\omega)\) for all \(n>N(\tilde\omega)\),
    and thus \(\argmax_i \mu^n_i(\BFx;\tilde\omega)\to i^*(\BFx;\tilde\omega)\) as \(n\to\infty\).
\end{proof}

\subsection{F. Proof of Theorem \ref{theo:consistency-cost-quasi}}

The steps to prove Theorem \ref{theo:consistency-cost-quasi} are exactly the same to those for \Cref{theo:consistency},
and we only need to modify the arguments related to the actual sampling decision (which is $(a^n,\BFv^n)$ satisfies \cref{eq:max-int-KG} in the IKG policy, and $(\tilde{a}^n,\tilde{\BFv}^n)$ satisfies \cref{eq:max-int-KG-quasi} in the quasi-IKG policy).
So, we will not repeat the entire proofs, but only point out the modification briefly.
Specifically, the two main steps to prove Theorem \ref{theo:consistency-cost-quasi} are summarized as the following \Cref{prop:step1-quasi,prop:step2-quasi}, which are parallel to \Cref{prop:step1,prop:step2}.

\begin{proposition}\label{prop:step1-quasi}
    Fix \(i=1,\ldots,M\). If  \Cref{assump:compact,assump:stationary,assump:var} hold and \(\eta_i^\infty=\infty\) a.s., then for any \(\BFx\in\calX\), \(k_i^\infty(\BFx,\BFx)=0\) a.s. under the quasi-IKG policy.
\end{proposition}

\begin{proof}[Proof of \Cref{prop:step1-quasi}]
All the intermediate lemmas for \Cref{prop:step1} directly apply to \Cref{prop:step1-quasi}, except for \Cref{lemma:liminf}.
Instead, we now need to show that under the quasi-IKG policy, for any \(\BFx\in\calX\),
$\liminf_{n\to\infty} \IKG^n(i,\BFx) = 0.$
We first observe that, for the sequence \(\{\tilde{\BFv}^n\in\calX:\tilde{a}^n=i,n=0,1,\ldots,\}\), we can still show \(\liminf_{n\to\infty} \IKG^{n}(i,\tilde{\BFv}^n)  = 0\) with the same arguments.
Then, by the definition of the quasi-IKG policy \eqref{eq:max-int-KG-quasi},
$$\IKG^{n}(i,\tilde{\BFv}^n) = \IKG^{n}(\tilde{a}^n,\tilde{\BFv}^n) \geq \IKG^n(a^n,\BFv^n) - \varepsilon_n \geq \IKG^n(i,\BFv^n) - \varepsilon_n .$$
On the other hand, since \(\IKG^{n}(i,\BFv^n) = \max_{\BFx\in\calX}\IKG^{n}(i,\BFx)\),
then for any \(\BFx\in\calX\), \cref{eq:liminfIKG} is replaced by
$$
0 \leq \liminf_{n\to\infty} \IKG^n(i,\BFx) \leq\liminf_{n\to\infty} \IKG^{n}(i,\BFv^n) \leq  \liminf_{n\to\infty} [\IKG^{n}(i,\tilde{\BFv}^n) + \varepsilon_n] = 0,
$$
where the equality is due to \(\liminf_{n\to\infty} \IKG^{n}(i,\tilde{\BFv}^n)  = 0\) and $\varepsilon_n \to 0$.
Finally, the proof of \Cref{prop:step1-quasi} follows the similar arguments as in the proof of \Cref{prop:step1}.
\end{proof}

\begin{proposition}\label{prop:step2-quasi}
    If  \Cref{assump:compact,assump:stationary,assump:var} hold, then \(\eta_i^\infty=\infty\) a.s.  for each \(i=1,\ldots,M\) under the quasi-IKG policy.
\end{proposition}

\begin{proof}[Proof of \Cref{prop:step2-quasi}]
All the intermediate lemmas for \Cref{prop:step2} directly apply to \Cref{prop:step2-quasi}.
We then proceed by following the same arguments as in the proof \Cref{prop:step2}, with $T_i(\omega_0)$ meaning that the quasi-IKG policy does not choose alternative \(i\) for \(n > T_i(\omega_0)\).
After we obtain \cref{eq:contradiction_A}, we now want to show that quasi-IKG policy must choose alternative \(i\in A\) at some time \(\tilde n(\omega_0) > T(\omega_0)\), which leads to the contradiction.

Due to \cref{eq:contradiction_A}, there exists some small $\Delta>0$ such that
    $$
    \min_{i\in A} \max_{\BFx\in\calX} \ [c_i(\BFx)]^{-1} Q^\dag(S^\infty(\omega_0), i,\BFx) - \Delta
        > \max_{i'\notin A} \max_{\BFx\in\calX} \ [c_{i'}(\BFx)]^{-1} Q^\dag(S^\infty(\omega_0), i',\BFx)=0.
    $$
Notice that \(S^n(\omega_0) \to S^\infty(\omega_0)\) pointwise as \(n\to\infty\) since \(\omega_0\in\Omega_1\).
Hence, there exists a finite number \(\tilde n_1(\omega_0) > T(\omega_0)\) such that
\begin{equation}\label{eq:Sn-large-n}
\min_{i\in A}\max_{\BFx\in\calX} \ [c_i(\BFx)]^{-1} Q^\dag(S^n, i,\BFx) - \Delta/2
        > \max_{i'\notin A} \max_{\BFx\in\calX} \ [c_{i'}(\BFx)]^{-1} Q^\dag(S^n, i',\BFx),
\end{equation}
for all $n\geq \tilde n_1(\omega_0)$.
Since $\varepsilon_n\to 0$ as \(n\to\infty\), there exists a finite number \(\tilde n_2(\omega_0) > T(\omega_0)\) such that $\varepsilon_n < \Delta/2$ for all $n\geq \tilde n_2(\omega_0)$.
Then, we can conclude that at time \(\tilde n(\omega_0) \coloneqq \max(\tilde n_1(\omega_0),\tilde n_2(\omega_0)) > T(\omega_0)\), quasi-IKG policy must choose alternative \(i\in A\).
Otherwise, for $n=\tilde n(\omega_0)$, if $\tilde{a}^n \notin A$, then by \cref{eq:IKG-3,eq:Sn-large-n},
$$
\IKG^n(\tilde{a}^n,\tilde{\BFv}^n) < \IKG^n(a^n,\BFv^n)  - \Delta/2 < \IKG^n(a^n,\BFv^n) -\varepsilon_n,
$$
which violates the definition of quasi-IKG policy defined in \cref{eq:max-int-KG-quasi}.
\end{proof}

\subsection{G. Gradient Calculation}

It is easy to see that
\begin{align*} 
    g_i^n(\BFv, \BFx)                                & = \pdv{[h_i^n(\BFv, \BFx)/c_i(\BFx)]}{\BFx} = \frac{\pdv{h_i^n(\BFv, \BFx)}{\BFx} c_i(\BFx) - h_i^n(\BFv, \BFx) \dv{c_i(\BFx)}{\BFx}}{[c_i(\BFx)]^2}, \\[0.5ex]
    \pdv{h_i^n(\BFv, \BFx)}{\BFx} & =
    \begin{cases}
        \phi \qty(\abs{\frac{\Delta_i^n(\BFv)}{\tilde\sigma_i^n(\BFv,\BFx)} }) \pdv{\tilde\sigma_i^n(\BFv,\BFx)}{\BFx},  & \text{if } \tilde\sigma_i^n(\BFv,\BFx) > 0, \\
        0,                                                                                                                                                & \text{if } \tilde\sigma_i^n(\BFv,\BFx) = 0, \\
        -\phi \qty(\abs{\frac{\Delta_i^n(\BFv)}{\tilde\sigma_i^n(\BFv,\BFx)} }) \pdv{\tilde\sigma_i^n(\BFv,\BFx)}{\BFx}, & \text{if } \tilde\sigma_i^n(\BFv,\BFx) < 0,
    \end{cases}
\end{align*}
and that, by the definition of \(\tilde\sigma_i^n(\BFv,\BFx)\) in \cref{eq:sigma} and \Cref{assump:stationary,assump:var},
\begin{equation}\label{eq:partial_sigma_0}
    \pdv{\tilde\sigma_i^n(\BFv,\BFx)}{\BFx} = [k_i^n(\BFx,\BFx)+\lambda_i(\BFx)]^{-\frac{1}{2}}\pdv{k^n_i(\BFv,\BFx)}{\BFx} - \frac{[k_i^n(\BFx,\BFx)+\lambda_i(\BFx)]^{-\frac{3}{2}}k^n_i(\BFv,\BFx)}{2} \left[\dv{k^n_i(\BFx,\BFx)}{\BFx} + \dv{\lambda_i(\BFx)}{\BFx} \right],
\end{equation}
provided that the prior correlation function \(\rho_i\) and the cost function \(c_i\) are both differentiable.
Assuming \(\rho_i\) to be differentiable excludes some covariance functions that satisfy \Cref{assump:stationary} such as the \(\text{Mat\'{e}rn}(\nu)\) type with \(\nu=1/2\), but many others including both the \(\text{Mat\'{e}rn}(\nu)\) type with \(\nu>1\) and the SE type do have the desired differentiability. We next calculate analytically the derivatives of \(k_i^n(\cdot,\cdot)\) in \cref{eq:partial_sigma_0} for several common covariance functions. The calculation is  a routine exercise so we omit the details.

Throughout the subsequent \Cref{ex:SE,ex:Matern,ex:Matern2}, we use the following notation. For \(i=1,\ldots,M\),
\[\BFalpha_i \coloneqq (\alpha_{i,1},\ldots,\alpha_{i,d})^\intercal\qq{and} r_i(\BFx,\BFx') \coloneqq \sqrt{\sum_{j=1}^d\alpha_{i,j} (x_j-x_j')^2 }.\]
Recall that \(\BFV_i^n\) denotes the set of locations of the samples taken from \(\theta_i\) up to time \(n\).
With slight abuse of notation, here we treat \(\BFV_i^n\) as a matrix wherein the columns are corresponding to the points in the set and arranged in the order of appearance.
Moreover, for notational simplicity, let \(\BFV_i^n \coloneqq (\BFv_1,\ldots,\BFv_{m_i^n})\), where \(m_i^n\) is the number of columns of \(\BFV_i^n\).
Let \(\BFX\) be a matrix with the same dimension as \(\BFV_i^n\) and all columns are identically \(\BFx\).
We adopt the denominator layout for matrix calculus.
For the following \Cref{ex:SE,ex:Matern,ex:Matern2}, it can be shown that
\begin{align}
    \pdv{k_i^n(\BFv,\BFx)}{\BFx} & = \mathrm{diag}(\BFalpha_i) (\BFx - \BFv) a_0 - \BFA k_i^0(\BFV_i^n,\BFv), \label{eq:dkvx} \\
    \dv{k_i^n(\BFx,\BFx)}{\BFx}           & = -2 \BFA k_i^0(\BFV_i^n,\BFx), \label{eq:dkxx}
\end{align}
where
\[\BFA \coloneqq \mathrm{diag}(\BFalpha_i) (\BFX - \BFV_i^n) \mathrm{diag}\{a_1,\ldots,a_{m_i^n}\} [k_i^0(\BFV_i^n,\BFV_i^n)+\lambda_i\BFI]^{-1},\]
while the values of \(a_0, a_1, \ldots, a_{m_i^n}\) depend on the choice of the covariance function.

\begin{example}[SE]\label{ex:SE}
    Let  \(k_i^0(\BFx,\BFx')=\tau_i^2\exp\qty(-r_i^2(\BFx,\BFx'))\).
    Then, in \cref{eq:dkvx,eq:dkxx}, \(a_0 \coloneqq -2 k_i^0(\BFv,\BFx)\) and \(a_\ell \coloneqq -2 k_i^0(\BFv_\ell,\BFx)\), for \(\ell = 1,\ldots,m_i^n\).
\end{example}

\begin{example}[\(\text{Mat\'{e}rn}(3/2)\)]\label{ex:Matern}
    Let \(k_i^0(\BFx,\BFx') = \tau_i^2 \qty(1+\sqrt{3}r_i(\BFx,\BFx'))\exp \qty( - \sqrt{3}r_i(\BFx,\BFx'))\).
    Then, in \cref{eq:dkvx,eq:dkxx},
    \[ a_\ell \coloneqq \sqrt{3} r_i^{-1}(\BFv_\ell,\BFx) \qty[\tau_i^2 \exp \qty(-\sqrt{3}r_i(\BFv_\ell,\BFx)) - k_i^0(\BFv_\ell,\BFx) ], \]
    for \(\ell = 0,1,\ldots,m_i^n\) and \(\BFv_0 = \BFv\).
\end{example}

\begin{example}[\(\text{Mat\'{e}rn}(5/2)\)]\label{ex:Matern2}
    Let \(k_i^0(\BFx,\BFx') = \tau_i^2 \qty(1+\sqrt{5}r_i(\BFx,\BFx')+\frac{5}{3}r_i^2(\BFx,\BFx'))\exp \qty(- \sqrt{5}r_i(\BFx,\BFx'))\). Then, in \cref{eq:dkvx,eq:dkxx},
    \[ a_\ell \coloneqq \qty(\sqrt{5} r_i^{-1}(\BFv_\ell,\BFx) + {\textstyle \frac{10}{3}}) \tau_i^2 \exp \qty(-\sqrt{5}r_i(\BFv_\ell,\BFx)) - \sqrt{5} r_i^{-1}(\BFv_\ell,\BFx)k_i^0(\BFv_\ell,\BFx), \]
    for \(\ell = 0,1,\ldots,m_i^n\) and \(\BFv_0 = \BFv\).
\end{example}

\subsection{H. Implementation Issues}

A plain-vanilla implementation of \(\widehat{\IKG}^n(i,\BFx)\) in \cref{eq:sample_average} may encounter rounding errors, since \(h_i^n(\BFxi_j,\BFx)\) may be rounded to zero when evaluated via \cref{eq:function_h_i}; see~\cite{FrazierPowellDayanik09} for discussion on a similar issue.
To enhance numerical stability, we first evaluate the logarithm of the summand and then do exponentiation.
For notational simplicity, we set
\[u_j \coloneqq \abs{\Delta_i^n(\BFxi_j)/\tilde\sigma_i^n(\BFxi_j, \BFx)} \qq{and} h_i^n(\BFxi_j, \BFx) = \abs{\tilde\sigma_i^n(\BFxi_j, \BFx)} [\phi(u_j) - u_j \Phi(-u_j)].\]
If \(\tilde\sigma_i^n(\BFxi_j, \BFx)\neq 0\),  we compute
\[
    g_j \coloneqq \log\qty(\frac{h_i^n(\BFxi_j, \BFx)}{J} ) = \log\qty(\frac{\abs{\tilde\sigma_i^n(\BFxi_j, \BFx)}}{\sqrt{2\pi} J}) - \frac{1}{2} u_j^2 + \log\qty(1 - u_j \frac{\Phi(-u_j)}{\phi(u_j)}),
\]
where \({\Phi(-u_j)}/{\phi(u_j)}\) is known as the Mills ratio, and can be asymptotically approximated by \(u_j/(u_j^2+1)\) for large \(u_j\). Moreover, \( \log \qty(1 + x )\) can be accurately computed by \texttt{log1p} function available in most numerical software packages.
At last, we compute
\[\log \widehat{\IKG}^n(i,\BFx) =  \log \sum_{j \in \calJ} e^{g_j} -\log(c_i(\BFx)) = g^* + \log \sum_{j \in \calJ} e^{g_j - g^*} -\log(c_i(\BFx)),\]
where \( \calJ \coloneqq \{j:\tilde\sigma_i^n(\BFxi_j, \BFx)\neq  0,\,j=1,\ldots,J\}\), and \(g^* = \max_{j \in \calJ} g_j\);
we set \(\log \widehat{\IKG}_{\gamma}^n(i,\BFx) = -\infty\) if \(\calJ\) is empty. The above procedure is summarized in \Cref{alg:ComputingIKG}.

\begin{algorithm}[ht]
    \caption{Computing \(\log \widehat{\IKG}^n(i,\BFx)\).}\label{alg:ComputingIKG}
    {\fontsize{9}{8}\selectfont
        \begin{algorithmic}[1]
            \Statex \hspace{-17pt} \textbf{Inputs:}
            \(\mu_1^n,\ldots,\mu_M^n, k_1^n,\ldots,k_M^n, \lambda_1,\ldots,\lambda_M, \BFxi_1, \ldots, \BFxi_J, i,\BFx, c_i(\BFx)\)
            \vspace{1pt}
            \Statex \hspace{-17pt} \textbf{Outputs:}
            \emph{log\_IKG}

            \State \(\calJ \gets \emptyset\), \emph{log\_IKG} \(\gets -\infty\)
            \For{\(j=1\) to \(J\)}
            \If{\(\abs{\tilde\sigma_i^n(\BFxi_j, \BFx)} > 0\)}
            \vspace{1pt}
            \State \(u \gets \abs{\Delta_i^n(\BFxi_j)} / \abs{\tilde\sigma_i^n(\BFxi_j, \BFx)}\)
            \vspace{1pt}
            \If{\(u < 20\)}
            \State \(r \gets {\Phi(-u)}/{\phi(u)}\)
            \Else
            \State \(r \gets u/(u^2+1)\)
            \EndIf
            \State \(g_j \gets \log \qty(\frac{\abs{\tilde\sigma_i^n(\BFxi_j, \BFx)}}{\sqrt{2\pi} J}) - \frac{1}{2} u^2 + \texttt{log1p}(- u r)\)
            \Comment{\(\texttt{log1p}(x) = \log(1+x)\).}
            \State \(\calJ \gets \{\calJ, j\}\)
            \EndIf
            \EndFor
            \vspace{1pt}
            \If{\(\calJ \neq \emptyset\)}
            \State \(g^* \gets \max_{j \in \calJ} g_j\)
            \State \emph{log\_IKG} \(\gets g^* + \log \sum_{j \in \calJ} e^{g_j-g^*} -\log(c_i(\BFx))\)
            \EndIf
        \end{algorithmic}
    }
\end{algorithm}

In the implementation of SGA, we adopt two well-known modifications.
\begin{enumerate}[label=(\roman*)]
    \item We use \emph{mini-batch} SGA to have more productive iterations.
    Specifically, in each iteration \cref{eq:SGA_ite}, instead of using a single \(g_i^n(\BFxi_k, \BFx_k)\) as the gradient estimate, we use the average of \(m\) independent estimates \(g_i^n(\BFxi_{k1}, \BFx_k), \ldots, g_i^n(\BFxi_{km}, \BFx_k)\), which is denoted as \(\bar{g}_i^n(\BFxi_{k1},\ldots,\BFxi_{km}, \BFx_k)\).
    \item We adopt the \emph{Polyak-Ruppert averaging}~\citep{PolyakJuditsky92} to mitigate of the algorithm's sensitivity on the choice of the step size.
    Specifically, when \(K\) iterations are completed, we report \(\frac{1}{K+2-K_0}\sum_{k=K_0}^{K+1} \BFx_k\), instead of \(\BFx_{K+1}\), as the approximated solution of \(\BFv^n_i\), where \(1 \leq K_0 \leq K\) is a pre-specified integer.
\end{enumerate}

\begin{algorithm}[ht]
    \caption{Approximately Computing \((a^n,\BFv^n)\) Using SGA.}\label{alg:SGA}
    {\fontsize{9}{8}\selectfont
        \begin{algorithmic}[1]
            \Statex \hspace{-17pt} \textbf{Inputs:}
            \(\mu_1^n,\ldots,\mu_M^n, k_1^n,\ldots,k_M^n, \lambda_1,\ldots,\lambda_M, \BFxi_1, \ldots, \BFxi_J, c_1,\ldots,c_M\)
            \vspace{1pt}
            \Statex \hspace{-17pt} \textbf{Outputs:}
            \(\hat{a}^n,\hat{\BFv}^n\)
            \For{\(i=1\) to \(M\)}
            \State \(\BFx_1 \gets\) initial value
            \For{\(k=1\) to \(K\)}
            \State Generate independent sample \(\{\BFxi_{k1},\ldots,\BFxi_{km}\}\) from density \(\gamma(\cdot)\)
            \State \(\BFx_{k+1} \gets \Pi_{\calX} \qty[\BFx_k + b_k \bar{g}_i^n(\BFxi_{k1},\ldots,\BFxi_{km}, \BFx_k) ]\)
            \Comment{Mini-batch SGA.}
            \EndFor
            \State \(\hat{\BFv}^n_i \gets \frac{1}{K+2-K_0}\sum_{k=K_0}^{K+1} \BFx_k\)
            \Comment{Polyak-Ruppert averaging.}
            \State \emph{log\_IKG}\(_i\) \(\gets \log \widehat{\IKG}^n(i,\hat{\BFv}^n_i)\)
            \Comment{Call \Cref{alg:ComputingIKG}.}
            \EndFor
            \State \(\hat{a}^n \gets \argmax_i\) \emph{log\_IKG}\(_i\)
            \State \(\hat{\BFv}^n \gets \hat\BFv^n_{\hat{a}^n}\)
        \end{algorithmic}
    }
\end{algorithm}

Upon computing \(\hat\BFv_i^n\approx\argmax_{\BFx}\IKG^n(i,\BFx)\) with SGA for each \(i\), we set
    \[\hat a^n= \argmax_{1\leq i\leq M} \log \widehat{\IKG}^n(i,\hat\BFv_i^n)\qq{and} \hat \BFv^n= \hat\BFv_{\hat a^n}^n,\]
to be the sampling decision at time \(n\), i.e., let \((\hat a^n, \hat \BFv^n)\) be the computed solution for \((a^n,\BFv^n)\) under the IKG policy. The complete procedure is summarized in \Cref{alg:SGA}.

\subsection{I. Additional Numerical Experiments}

\subsubsection{Computational cost comparison}

We conduct simple experiment to compare the computational cost when the IKG sampling policy defined in \cref{eq:max-int-KG}, i.e., $(a^n,\BFv^n) \in \argmax_{1\leq i\leq M, \BFx\in\calX}\IKG^n(i,\BFx)$, is solved purely using sample average approximation (SAA) method or our proposed SGA (together with SAA).
For IKG with SGA, here refer to as method 1, as described in \Cref{sec:numerical}, the computation of \cref{eq:max-int-KG} consists of two steps.
Step (i) is to solve \(\BFv^{n}_i = \max_{\BFx\in\calX} \IKG^n(i,\BFx)\) for all \(i=1,\ldots,M\) with SGA,
and step (ii) is to solve \(a^n=\argmax_{1\leq i\leq M} \IKG^n(i, \BFv^{n}_i)\) with SAA.
For IKG with pure SAA, here refer to as method 2, the problems in the above two steps are both solved with SAA.
In particular, the problem in step (i) is converted into a continuous deterministic optimization after applying SAA, which is solved directly using the \texttt{fmincon} solver in MATLAB.
It is expected that for either method, the computational cost will increase as the dimensionality $d$ increases.
But to evaluate the exact value of the computational cost, one needs to know the true optimal solution of \cref{eq:max-int-KG} and control the optimality gap when a specific method is used.
Here we simply compare the relative computational cost of two methods by roughly controlling the resulting OC at the same level.

The same problem in \Cref{sec:numerical} is considered and we let the dimensionality $d$ increase from 1 to 7.
The density function of covariates is the uniform distribution and the cost function is constantly 1.
All the parameters for the problem, the IKG with SGA (i.e., method 1) and the evaluation of OC are the same as before.
For each $d$, the sampling policy under the two methods is carried out respectively until the budget $B=100$ is exhausted, and the $\widehat{\text{OC}}(B)$ curves are obtained.
For fair comparison, we let the step (ii) of method 2 be exactly the same (i.e., same sample used) as that of method 1, and tune the sample size in step (i) of method 2 as follows.
We gradually increase the sample size of random covariates used in SAA (not the number of samples from alternatives), until the resulting $\widehat{\text{OC}}(B)$ curve from method 2 is roughly comparable to that from method 1.
Note that for either method the computation time needs to solve \cref{eq:max-int-KG} depends on the number of samples allocated to each alternative so far, and will increase as the samples accumulate.
So, we report the total computation time spent on solving \cref{eq:max-int-KG} during the entire sampling process (until the budget $B=100$ is exhausted, which means \cref{eq:max-int-KG} is solved for 100 times), averaged on $L=30$ replications, for the two methods, which are denoted as $T_1(d)$ and $T_2(d)$ respectively.

\begin{figure}[!htp]
\centering
\caption{Computational cost comparison between SGA and SAA.} \label{fig:computation_cost}
\includegraphics[width=0.40\textwidth]{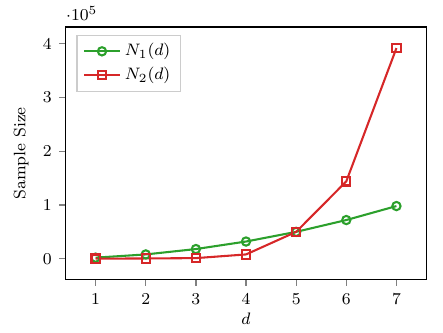}
\hspace{0.05\textwidth}
\includegraphics[width=0.40\textwidth]{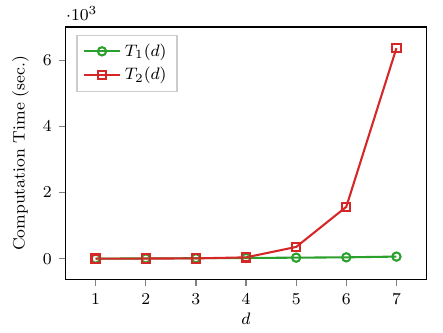}
\end{figure}

The following \Cref{fig:computation_cost} shows the comparison between methods 1 and 2.
Left panel of \Cref{fig:computation_cost} shows the sample sizes of random covariates used to approximately solve \(\BFv^{n}_i = \max_{\BFx\in\calX} \IKG^n(i,\BFx)\) for each $i$ in step (i) by the two methods, which are denoted as $N_1(d)$ and $N_2(d)$ respectively.
Note that $N_1(d)= mK =20d \times 100d = 2000d^2$ as specified, and $N_2(d)$ is tuned so that the performance of method 2 matches that of method 1.
Right panel of \Cref{fig:computation_cost} shows the computation times $T_1(d)$ and $T_2(d)$ (in MATLAB, Windows 10 OS, 3.60 GHz CPU, 16 GB RAM).
It can be seen that IKG with SGA (i.e., method 1) scales much better in dimensionality $d$ than IKG with pure SAA (i.e., method 2).
Recall that for each method, the sample size of covariates in step (ii) is set as $500d^2$.
So, even consider the fact that in method 2 the function approximation in step (i) can be directly used in step (ii), which saves the sample of covariates and the relevant computation in step (ii), the entire sample size and the computation time of method 2 still grows much faster than method 1.

\subsubsection{Estimated sampling variance}

In practice, the sampling variance \(\lambda_i(\BFx)\) is usually unknown and needs to be estimated.
We suggest to follow the approach in \cite{AnkenmanNelsonStaum10}.
Specifically, for each alternative $i$, at some predetermined design points $\BFx^1,\ldots,\BFx^m$, multiple simulations are run and the sample variances are computed, which are denoted as $s_i^2(\BFx^1),\ldots,s_i^2(\BFx^m)$.
Then ordinary kriging (i.e., Gaussian process interpolation) is used to fit the entire surface of \(\lambda_i(\BFx)\).
Under the Bayesian viewpoint, it is equivalent to impose a Gaussian process with constant mean function \(\mu_i^0(\BFx)\equiv\mu_i^0\) and covariance function \(k_i^0(\BFx,\BFx')\) as prior of \(\lambda_i(\BFx)\), and compute the posterior mean function by ignoring the sampling variance at the design points, i.e.,
$$
\mu_i^m(\BFx) = \mu_i^0 +  k_i^0(\BFx,\BFX_i)k_i^0(\BFX_i,\BFX_i)^{-1}[\BFy_i-\mu_i^0 \BFI],
$$
where \(\BFX_i \coloneqq (\BFx^1,\ldots,\BFx^m)\) and $\BFy_i \coloneqq (s_i^2(\BFx^1),\ldots,s_i^2(\BFx^m))^\intercal$.
Then, $\mu_i^m(\BFx)$ is used as estimate of \(\lambda_i(\BFx)\), and the IKG policy is applied as if \(\lambda_i(\BFx)\) was known.
In ordinary kriging, $\mu_i^0$ and the parameters in \(k_i^0(\BFx,\BFx')\) are usually optimized via maximum likelihood estimation (MLE).

\begin{figure}[!h]
    \begin{center}
        \caption{Estimated opportunity cost (vertical axis) as a function of the sampling budget (horizontal axis).} \label{fig:estimate_var}
            \includegraphics[width=0.36\textwidth]{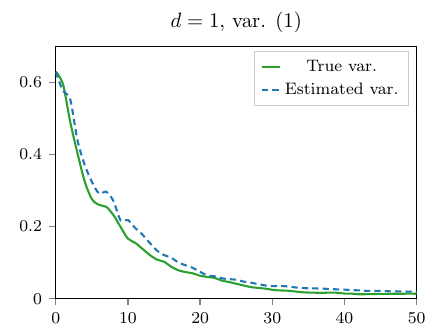} \hspace{0.05\textwidth}
            \includegraphics[width=0.36\textwidth]{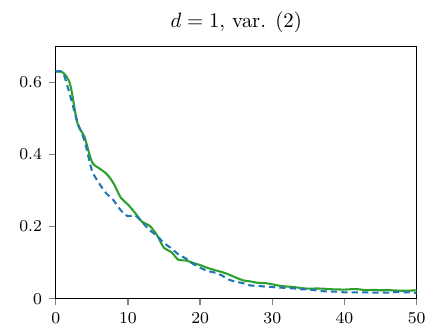}
            \includegraphics[width=0.36\textwidth]{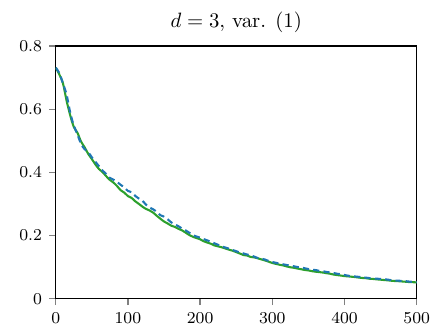} \hspace{0.05\textwidth}
            \includegraphics[width=0.36\textwidth]{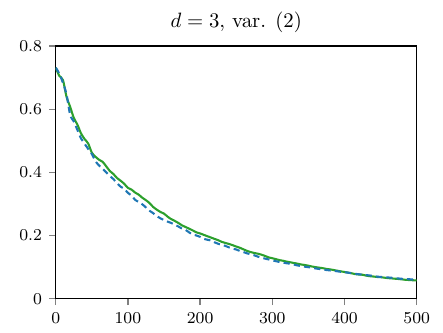}
    \end{center}
\end{figure}

We again consider the problem in \Cref{sec:numerical}.
To better investigate the effect of estimating \(\lambda_i(\BFx)\), we consider two sampling variance:
(1) \(\lambda_i(\BFx)\equiv 0.01\), as before;
(2) \(\lambda_i(\BFx) = 0.01 \times (1.5^{d-1} + \theta_i(\BFx))\).
The density function of covariates is the uniform distribution and the cost function is constantly 1.
All the other parameters for the problem are the same as before.
The prior for estimating \(\lambda_i(\BFx)\) is directly set as \(\mu_i^0= 0\), and \(k_i^0(\BFx,\BFx') = \exp\qty(- \frac{1}{d} \norm{\BFx-\BFx'}^2)\), without invoking the MLE.
The design points are generated by Latin hypercube sampling and the same design points are used for each $i$.
All the parameters for the IKG policy and the evaluation of OC are the same as before.
The following \Cref{fig:estimate_var} shows the estimated opportunity cost when the sampling variance is known or estimated using the above approach, for the case of $d=1$ or $3$ and sampling variance (1) or (2).
Numerical results show that the effect of estimating the sampling variance is minor for this problem, which agrees with the observation in \cite{AnkenmanNelsonStaum10}.

